\documentclass[12pt]{article}

\usepackage{amsmath,a4wide,amssymb,amsthm}
\usepackage{graphicx}

\newtheorem{theorem}{Theorem}
\newtheorem{lemma}[theorem]{Lemma}
\newtheorem{prop}[theorem]{Proposition}

\newtheorem{defi}[theorem]{Definition}
\newtheorem{ass}[theorem]{Assumption}

\numberwithin{equation}{section}
\numberwithin{theorem}{section}

\newcommand{\bs}{\boldsymbol}

\title{\bf On $q$-functional equations\\
and excursion moments
}

\author{\sc Christoph~Richard\\
\\
{\it Fakult\"at f\"ur Mathematik, Universit\"at Bielefeld,}\\
{\it Postfach 10 01 31, 33501 Bielefeld, Germany}}        

\begin{document}

\maketitle

\begin{abstract}
We analyse $q$-functional equations arising from tree-like combinatorial 
structures, which are counted by size, internal path length, and certain 
generalisations thereof. The corresponding counting parameters are labelled
by a positive integer $k$. We show the existence of a joint limit distribution 
for these parameters in the limit of infinite size, if the size generating 
function has a square root as dominant singularity. The limit distribution 
coincides with that of integrals of $k$-th powers of the standard Brownian 
excursion. Our approach yields a recursion for the moments of the limit
distribution. It can be used to analyse asymptotic expansions of the 
moments, and it admits an extension to other types of singularity.
\end{abstract}

Keywords: simply generated trees, $q$-difference equation, Brownian excursion, 
limit distribution

\section{Introduction and main results}

\subsection{A combinatorial motivation}

When studying combinatorial classes, a functional equation of the form
\begin{equation}\label{form:geneq}
G(x)=P(x,G(x))
\end{equation}
frequently arises, where $G(x)$ is the generating function of the
class, and $P(x,y)$ is a formal power series in two variables with real
coefficients.
Prominent examples are classes of simply generated trees, counted by 
number of vertices \cite{MM78}, classes of directed square lattice paths 
counted by length \cite{BF02}, and classes of square
lattice polygons, counted by perimeter \cite{Bou96}. Indeed, there exist 
combinatorial bijections between corresponding models of trees, paths, 
and polygons, see e.g.~\cite{S99} for Catalan trees, Dyck paths, 
and staircase polygons. Eqn.~\eqref{form:geneq} reflects a combinatorial 
decomposition of the given class. If $P(0,y)\equiv0$ and $\frac{\partial P}{\partial y}
(0,y)\not\equiv1$, then eqn.~\eqref{form:geneq} admits a unique solution 
$G(x)\in\mathbb R[[x]]$ such that $G(0)=0$. Here, $\mathbb R[[x]]$ denotes the 
ring of formal power series in $x$ with coefficients in $\mathbb R$. Often, the 
coefficients of $P(x,y)$ are non-negative real numbers. Then, the series $G(x)$ is a power 
series with non-negative coefficients, typically analytic at $x=0$ with a
square root as dominant singularity, see e.g.~\cite[Thm.~10.6]{O95}, \cite{D97},
or \cite[Ch.~VII.4]{FS07}, and references therein.

\smallskip

We are interested in certain deformations of the above equation. This is done by 
introducing a new variable $q$, such that the limit $q\to1$ reduces to
the original equation. For example, the functional equation 
\begin{equation}\label{form:geneq2}
G(x,q)=P(x,G(qx,q))
\end{equation}
defines a formal power series in $x$ with polynomial coefficients, i.e., 
$G(x,q)\in\mathbb R[q][[x]]$, if eqn.~\eqref{form:geneq} defines a power series 
$G(x)\in\mathbb R[[x]]$. The above equation is a \emph{$q$-difference equation},
see e.g.~\cite{VRSZ03} and references therein. It appears 
in classes of simply generated trees, counted by number of vertices and internal 
path length (i.e., the sum of the vertex distances to the root), 
in classes of directed square lattice paths, counted 
by length and area under the path, and in classes of square
lattice polygons, counted by perimeter and area. For some models, 
an explicit expression for its generating function $G(x,q)$ is known, see 
e.g.~\cite{Bou96,PB95a,PO95} and references therein. Such an expression 
typically contains $q$-products and has a natural boundary $|q|=1$. 
An interesting question concerns the statistics of the additional 
counting parameter, e.g., in a uniform ensemble in the limit of large system size. 
It is known (\cite{T91, D99}, see also \cite{R02}) that eqn.~\eqref{form:geneq2} leads, for certain 
$q$-difference equations and after appropriate normalisation, to the \emph{Airy 
distribution} as the limit distribution for the 
additional parameter. This distribution is known to also describe the area under a 
Brownian excursion, see the following subsection. Note that generalisations 
of eqn.~\eqref{form:geneq} other than eqn.~\eqref{form:geneq2} have also been 
studied previously. A class of equations, 
which leads to Gaussian limit laws, is discussed in \cite{D97}.

\smallskip

The idea of iterating the above deformations has been considered by Duchon 
\cite{D98,D99}. The deformation variables may be denoted by $q_k$, where 
$k\in\{1,\ldots,M\}$, and the associated counting parameters 
are called \emph{parameters of rank $k+1$}. For eqn.~\eqref{form:geneq2}, 
an example is given by
\begin{equation}\label{form:geneq3}
G(x,q_1,\ldots,q_M)=P(x,G(xq_1\cdot\ldots\cdot q_M,q_1q_2\cdot\ldots\cdot q_M,
q_2q_3\cdot\ldots\cdot q_M,\ldots,q_M)).
\end{equation}
Here, $G(x,q_1,\ldots, q_M)=\sum_n p_n(q_1,\ldots, q_M) x^n$ is a formal power
series in $x$, with polynomial coefficients $p_n(q_1,\ldots, q_M)\in\mathbb 
R[q_1,\ldots, q_M]$. The name ``parameter of rank $k+1$'' reflects that 
$k$ is the smallest integer $r$, such that the degree of the polynomial 
$p_n(q_1,\ldots, q_M)$ in $q_k$ is bounded by 
$c\,n^{r+1}$ for some constant $c$ (see \cite[Lemma~1]{D99} and \cite{D98}). 
We will call an equation like eqn.~\eqref{form:geneq3} a \emph{$q$-functional equation} 
(Definition \ref{defi:qfunc}). 
Again the question arises, under which conditions a limit distribution for the additional
counting parameters exists. 
In this paper, we shall show this for a class of deformations 
which we call \emph{$q$-shift} (Definition \ref{def:qshift}), the main assumption 
on the $q$-functional equation being that the solution of the undeformed equation, 
eqn.~\eqref{form:geneq}, is analytic at the origin, with a square root as dominant singularity 
(Assumption \ref{ass}). See Theorem~\ref{theo:probdist} for a precise statement. The 
resulting limit distributions appear to be related to distributions of integrals of $k$-th 
powers of the standard Brownian excursion. We will obtain a recursion for the moments 
of the joint distribution.  Our approach is based on the multivariate \emph{moment method}
(see e.g.~\cite{JLR00, Bill95}). The univariate case 
$M=1$ has been studied previously \cite{T91, FPV98, D99, R02}, and recursions for $M=2$ have been 
studied in \cite{NT03b}.

\smallskip

Before we consider general $q$-functional equations in Section~\ref{sec:qfunc}, let us 
first discuss the above questions in more detail for the simple example of Dyck paths.

\subsection{Dyck paths and Brownian excursions}\label{sec:Dyck}

We review the connection between Dyck paths and Brownian excursions.
This relates, in particular, the moments of height
of a random Dyck path to Brownian excursion moments.
Since the generating function of Dyck paths, counted by length and moments 
of their height, provides a simple example of a solution of a  $q$-functional 
equation, this also motivates the appearance of excursion moments in 
$q$-functional equations. We will state a central result of the paper in 
Theorem~\ref{theo:excmom}, which characterises the Brownian excursion random 
variables by a recursion for their moments.

\medskip

\begin{figure}[htb]
\begin{center}
\begin{minipage}[b]{0.65\textwidth}
\center{\includegraphics[width=8cm]{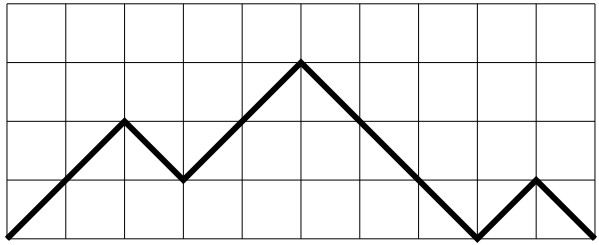}}
\end{minipage}
\end{center}
\caption{\label{fig:dyck}
\small A Dyck path of length $2n=10$. It is a sequence of two arches of lengths 8 and 2.}
\end{figure}

Let $\mathbb N_0=\mathbb N \cup \{0\}$ and $\mathbb R_{\ge0}=\{x\in\mathbb 
R: x\ge0\}$. A {\it Dyck path of length $2n$}, where $n\in \mathbb N_0$, is a map 
$y:[0,2n]\to\mathbb R_{\ge0}$, where $y(0)=y(2n)=0$ and $|y(i)-y(i+1)|=1$ 
for $i\in\mathbb N_0$ and $i<2n$. For non-integral argument, $y(s)$ is 
defined by linear interpolation. The values $y(s)$ are called the \emph{heights} 
of the path. An {\it arch of length $2n$} is a Dyck path of length $2n$, where $n>0$ and 
$y(s)>0$ for $s\in (0,2n)$. An example is given in Figure \ref{fig:dyck}.

\medskip

We are interested in a probabilistic interpretation of Dyck paths in a uniform
ensemble where, for fixed length, each of the finitely many Dyck paths occurs with 
the same probability. For $0\le s\le 2n$, let $\widetilde Y_{n}(s)$ denote the height 
of a random Dyck path of length $2n$. It is well known, see e.g.~\cite{BKR72}, that 
the average maximal height $\widetilde h(n)$ of a random Dyck path of length $2n$ has 
the asymptotic form $\widetilde h(n)\sim \sqrt{\pi n}$ as $n\to\infty$. 
In order to obtain a finite positive limit as $n$ approaches infinity, 
and in order to normalise the domain, we introduce the normalised height 
$Y_{n}(t)=(2n)^{-1/2}\widetilde Y_{n}(2nt)$, where $0\le t\le 1$. 

The sequence $\{Y_{n}(t)\}_{n\in\mathbb N}$ is a sequence of stochastic 
processes defined on $(C[0,1],{||\cdot ||_{\infty}})$ with the Borel 
$\sigma$-algebra, which converges in distribution to the standard Brownian 
excursion $e(t)$ of duration $1$, see \cite{Bel72,A93, MM03}. 
This implies convergence in distribution of sequences of continuous bounded 
functionals of $Y_n(t)$ towards the corresponding excursion functionals. 
(We refer to \cite{B99} for background about convergence of probability measures.)
There 
is a more recent result \cite[Thm.~9]{D03} that also
\begin{equation*}
\lim_{n\to\infty}\mathbb E_n[F(Y_{n}(t))] = \mathbb E[F(e(t))]
\end{equation*}
for continuous functionals $F:C[0,1]\to\mathbb R$ of \emph{polynomial} growth, 
i.e., for functionals such that there exists an $r\ge0$ with $|F(y)|\le 
||y||_\infty^r$ for all $y\in C[0,1]$. The above property is called \emph{polynomial 
convergence}. In particular, polynomial convergence implies
\emph{moment convergence} for such a functional, i.e., convergence of the
sequence of moments of order $l$, $\{\mathbb E_n[F(Y_{n}(t))^l]\}_{n\in
\mathbb N}$, for every $l\in\mathbb N_0$.  
For the functional of polynomial growth $F(y)=\int_0^1 y^k(t) {\rm d}t$, 
$k\in\mathbb N$, we call the random variable $\int_0^1e^k(t){\rm d}t$ 
the \emph{$k$-th excursion moment}.

\medskip

As counting parameters for Dyck paths of length $2n$, we consider the parameters 
\begin{equation}\label{form:dyckpara}
x_{k,n}=\sum_{i=0}^{2n} y^k(i) \qquad (k=1,\ldots,M).
\end{equation}
These are sums of $k$-th powers of heights, which we call \emph{$k$-th moments of 
height}. The parameter $x_{1,n}$ is the area under the Dyck path, the parameter 
$x_{2,n}$ is called the \emph{moment of inertia} of the Dyck path \cite{NT03b}.
Let $\widetilde X_{k,n}$ denote the $k$-th moment of height of a random Dyck
path of length $2n$. In terms of $\widetilde Y_{n}(s)$, 
the random variables  $\widetilde X_{k,n}$ are expressed as
\begin{equation*}
\widetilde X_{k,n}=\sum_{i=0}^{2n} \widetilde Y_{n}^k(i).
\end{equation*}
The scaling of the average height of a random Dyck path with its length 
suggests considering normalised random variables $X_{k,n}$, defined by
\begin{equation}\label{form:DyckRV}
(X_{1,n},X_{2,n},\ldots, X_{M,n})= \left(\frac{\widetilde X_{1,n}}{n^{(1+2)/2}},
\frac{\widetilde X_{2,n}}{n^{(2+2)/2}},\ldots, 
\frac{\widetilde X_{M,n}}{n^{(M+2)/2}}\right).
\end{equation}
Then, in terms of $Y_{n}(t)$, the normalised random variables 
$X_{k,n}$ are given by
\begin{equation}\label{form:Dyckran}
X_{k,n}=2^{(k+2)/2}\sum_{i=0}^{2n} Y^k_{n}\left(\frac{i}{2n}\right)\cdot 
\frac{1}{2n}.
\end{equation}
Since the above sum is a Riemann sum, one is led to expect convergence
in distribution and moment convergence of the sequence $\{X_{k,n}\}_{n\in\mathbb N}$ to 
$2^{(k+2)/2}\int_0^1 e^k(t){\rm d}t$, due to the convergence properties of 
$\{Y_{n}(t)\}_{n\in\mathbb N}$. This is the statement of the following proposition.

\begin{prop}\label{prop:reldyckexc}
For $k\in \{1,\ldots,M\}$, let $X_k=\int_{0}^1 e^k(t)\,{\rm d}t$ denote the $k$-th 
excursion moment. The sequence of Dyck path random variables eqn.~\eqref{form:DyckRV} 
converges to the normalised excursion moments $(2^{(1+2)/2}X_1,\ldots,2^{(M+2)/2}X_M)$ 
in distribution,
\begin{equation}\label{form:DyEx}
(X_{1,n},X_{2,n},\ldots, X_{M,n}) \stackrel{d}{\longrightarrow} 
(2^{(1+2)/2}X_1,2^{(2+2)/2}X_2,\ldots,2^{(M+2)/2}X_M) \qquad (n\to\infty).
\end{equation}
We also have moment convergence.
\end{prop}

\begin{proof}

For $k\in\{1,\ldots,M\}$, define $Y_{k,n}=\int_0^1 Y^k_n(t){\rm d}t$. As stated
above, convergence in distribution and moment convergence holds for the 
sequence $\left\{(Y_{1,n},\ldots,Y_{M,n})\right\}_{n\in\mathbb N}$, 
see \cite{Bel72,A93} and \cite[Thm.~9]{D03}. We argue that the sequences 
$\{(X_{1,n},\ldots,X_{M,n})\}_{n\in\mathbb N}$ and $\{(2^{(1+2)/2}Y_{1,n},
\ldots,2^{(M+2)/2}Y_{M,n})\}_{n\in\mathbb N}$ converge in distribution and for 
moments to the same limit.

Fix $k\in\{1,\ldots,M\}$ and consider, for a Dyck path $y$ of length $2n$ 
and for $m\in\{0,1,\ldots,2n-1\}$, the elementary estimate $(\max\{0,y(m)-1\})^k
\le \int_m^{m+1}y^k(s){\rm d}s\le(y(m)+1)^k$. Summing over $m$ yields, together 
with the binomial theorem, the estimate
\begin{equation}\label{form:estmom}
\sum_{l=0}^k \binom{k}{l}(-1)^l x_{l,n}\le\int_0^{2n} y^k(s){\rm d}s
\le \sum_{l=0}^k \binom{k}{l} x_{l,n}.
\end{equation}
This estimate translates to the distribution functions of the corresponding 
random variables. In terms of the random variables $X_{k,n}$ and 
$Y_{n}(t)$, we get
\begin{equation}\label{form:cdis}
\sum_{l=0}^k \binom{k}{l}(-1)^l \frac{X_{l,n}}{n^{(k-l)/2}}\stackrel{d}{\le}
2^{(k+2)/2}\int_0^{1} Y_{n}^k(t){\rm d}t \stackrel{d}{\le} \sum_{l=0}^k 
\binom{k}{l}\frac{X_{l,n}}{n^{(k-l)/2}},
\end{equation}
where the superscript $\stackrel{d}{}$ indicates that the relation is to be 
understood via distribution functions. Since the sequence of random variables 
$\{X_{l,n}/n^{(l-k)/2}\}_{n\in\mathbb N}$ converges in probability to zero, 
for $l=0,\ldots,k-1$, the random variables on the l.h.s.~and on the r.h.s.~of 
eqn.~\eqref{form:cdis} converge in distribution to the same limit. We conclude 
that the sequences $\{X_{k,n}\}_{n\in\mathbb N}$ and $\{2^{(k+2)/2}\int_0^1 Y_{n}^k (t)
{\rm d}t\}_{n\in\mathbb N}$ converge in distribution to the same limit. 
Moment convergence follows similarly from eqn.~\eqref{form:estmom}. 

Since the above estimates also hold jointly in $k\in\{1,\ldots,M\}$, the 
statement of the proposition follows.
\end{proof}

\noindent \textbf{Remark.} 
The above statement concerns Dyck paths in a uniform ensemble. The same 
argument as above leads to analogous results for more general ensembles of 
Dyck paths, which arise from the depth first process of simply generated trees 
\cite{A93, MM03}. Compare Theorem~\ref{theo:probdist} for a further generalisation.

\medskip

Due to the above result, Brownian excursion functionals may be
studied via their discrete Dyck path counterparts. Below, we will analyse the 
asymptotic behaviour of the moments of the random variables 
$(X_{1,n},\ldots,X_{M,n})$ eqn.~\eqref{form:DyckRV}, using the $q$-functional 
equation which the generating function of Dyck paths obeys. This implies a 
certain recursion for the moments of the joint distribution of the 
excursion moments $(X_1,\ldots,X_M)$.
We have the following result, which is obtained by combining the statements of 
Proposition~\ref{prop:reldyckexc} and Proposition~\ref{prop:dyckmom} below.
For $M$-dimensional vectors, we use the abbreviation $\bs k=(k_1,\ldots,k_M)$ 
and write $\bs0=(0,\ldots,0)$. For $i\in\{1,\ldots,M\}$, the unit vector $\bs e_i$ is 
defined by $(\bs e_i)_j=\delta_{i,j}$ for $j=1,\ldots,M$, and $\bs l\le\bs k$ for 
vectors $\bs l=(l_1,\ldots,l_M)$ and $\bs k=(k_1,\ldots,k_M)$ means that $l_i\le k_i$ 
for $i=1,\ldots,M$. A multivariate power series with complex coefficients is called 
\emph{entire}, if it converges for arbitrary complex arguments.

\begin{theorem}[Excursion moments]\label{theo:excmom}

For $k\in\{1,\ldots, M\}$, let $X_k=\int_{0}^1 e^k(t)\,{\rm d}t$ denote the $k$-th 
excursion moment. The moments of the joint distribution of $(X_1,\ldots,X_M)$ 
are given by 
\begin{equation}\label{eqn:exmomexpl}
\frac{\mathbb E[X_1^{k_1}\cdot\ldots\cdot X_M^{k_M}]}{k_1!\cdot\ldots\cdot k_M!}
=\frac{\sqrt{2\pi}}{\Gamma(\gamma_{k_1,\ldots,k_M})}2^{\gamma_{k_1,\ldots,k_M}}
\frac{f_{k_1,\ldots,k_M}}{2},
\end{equation}
where $\gamma_{k_1,\ldots,k_M}=\gamma_{\bs k}=-1/2+\sum_{i=1}^M(1+i/2)k_i$,
and where $\Gamma(z)$ is the Gamma function. 
For ${\bs k}\ne {\bf 0}$, the numbers $f_{k_1,\ldots,k_M}=f_{\bs k}$ are characterised 
by the recursion
\begin{equation}\label{form:exc2}
\gamma_{{\bs k}-{\bs e}_1}
f_{{\bs k}-{\bs e}_1}+\sum_{i=1}^{M-1}2(i+1)(k_i+1)
f_{{\bs k}-{\bs e}_{i+1}+{\bs e}_i}+
\sum_{{\bf 0}\le{\bs l}\le{\bs k}}f_{\bs l}
f_{{\bs k}-{\bs l}}=0,
\end{equation}
with boundary conditions $f_{\bf 0}=-4$ and $f_{\bs k}=0$, if $k_j<0$ for 
some $j\in\{1,\ldots, M\}$. 
The moment generating function $\mathbb E[\mathrm{e}^{t_1 X_1+\ldots+t_M X_M}]$ 
is entire. Hence, the joint distribution is uniquely defined by its moments. \qed
\end{theorem}

\noindent \textbf{Remarks.} \textit{i)} The theorem asserts that the numbers 
$f_{\bs k}$ can be recursively computed from eqn.~\eqref{form:exc2}. To 
see this, we remark that the above equation has $-f_{\bs k}=2f_{\bs k}
f_{\bs 0}/8$ as a term, and that all other numbers $f_{\bs l}$ in the equation 
satisfy $\bs l\prec \bs k$ for a suitable total order $\prec$, which is 
specified below in Definition~\ref{def:totord}. \\
\textit{ii)} In probability theory, Louchard's theorem \cite{L84} leads
to a characterisation of a certain Laplace transform of the moment generating 
function, for some excursion functionals including excursion moments. It 
is however difficult to extract moment values or moment recursions from 
these expressions.  For $M=1$, this has been done in \cite{L85}.
The moment generating function of the marginal distribution for $M=2$ has been
obtained in \cite[Thm.~2.4]{NT03b}. Moment values or recursions for $M>2$ 
have apparently not been derived.\\
\textit{iii)} 
We obtain Theorem~\ref{theo:excmom} by studying the corresponding Dyck
path functionals. Our discrete approach has been used previously \cite{NT03b}.
It led to a combinatorial derivation \cite[Thm.~3.1]{NT03b} of Louchard's formula 
for integrals over Brownian excursion polynomials. It also led to the values 
$\mathbb E[(X_1)^{k_1} (X_2)^{k_2}]$, see \cite[Table~2]{NT03b}.
For arbitrary $M$, our result is announced in \cite{R04}. 
[A factor 1/2 is missing on the r.h.s.~of eqn.~(6) in \cite{R04}].\\
\textit{iv)} Assuming only convergence in distribution of $(X_{1,n},\ldots,X_{M,n})$, 
the above result can be used to provide an alternative proof of moment
convergence.

\medskip

Let us briefly discuss the moments of the marginal distributions for $M\le4$.
For $M=1$, the first few coefficients $\mathbb E[(X_1)^k]$ are $1$, 
$\frac{1}{4}\sqrt{2\pi}$, $\frac{5}{12}$, $\frac{15}{128}\sqrt{2\pi}$, 
$\frac{221}{1008}$, $\frac{565}{8192}\sqrt{2\pi}$. These are related to the 
\emph{Airy distribution} \cite{FL01,J07,KMM07}, i.e., the distribution of area under the 
Brownian excursion $\sqrt{8} \, e(t)$. Explicit expressions are 
known for the moment generating function, for the density of the distribution,
and for the moments.

For $M=2$, the first few moments $\mathbb E[(X_2)^k]$
are $1$, $\frac{1}{2}$, $\frac{19}{60}$, $\frac{631}{2520}$, 
$\frac{1219}{5040}$, $\frac{92723}{332640}$, $\frac{1513891}{4036032}$. 
An explicit expression for the corresponding moment generating function 
\cite{NT03b} is
\begin{equation*}
\mathbb E[\mathrm{e}^{tX_2}]=\sum_{k=0}^\infty\frac{\mathbb E[(X_2)^k]}{k!} t^k=
\left( \frac{\sqrt{2t}}{\sin(\sqrt{2t})}\right)^{3/2}.
\end{equation*}
The first few moments $\mathbb E[(X_3)^k]$ are $1$, 
$\frac{3\sqrt{2\pi}}{16}$, $\frac{207}{560}$, $\frac{11907\sqrt{2\pi}}{65536}$, 
$\frac{88655283}{108908800}$, $\frac{1165359069\sqrt{2\pi}}{1476395008}$, and 
for the first few moments $\mathbb E[(X_4)^k]$ we have  $1$, $\frac{1}{2}$, 
$\frac{251}{840}$, $\frac{288751}{1201200}$, $\frac{19093793}{76236160}$, 
$\frac{105169404203}{3259095840000}$. It remains an open problem to find
explicit expressions for the moment generating functions.

\subsection{Dyck paths and $q$-functional equations}

We discuss the functional equation which the generating function of Dyck 
paths obeys, when counted by length and $k$-th moments of height,
$k=1,\ldots,M$. We then state in Theorem~\ref{theo:probdist} a convergence 
result for the limit distribution of counting parameters related to general 
$q$-functional equations. This is the main result of our paper.

\medskip

Let $\bs u=(u_0,u_1,\ldots,u_M)$ denote formal (commutative) variables. For a Dyck 
path $p$ of length $2n$, its \emph{weight} $w(p)$ is given by $w(p)=
u_0^{n_{}} u_1^{x_{1,n}}\cdot\ldots\cdot u_M^{x_{M,n}}$, where $x_{k,n}$
is the $k$-th moment of height eqn.~\eqref{form:dyckpara}. For the set $\cal D$ 
of Dyck paths, its generating function is the formal power 
series $D({\bs u})=\sum_{p\in\cal D} w(p)$. For the set $\cal A$ of 
arches, its generating function is the formal power series $A({\bs 
u})=\sum_{p\in\cal A} w(p)$. The variables $u_1,\ldots, u_M$ are interpreted 
as deformation variables. They may all be set equal to unity, in which case 
the generating function reduces to that of Dyck paths (arches), counted by length only.
The generating functions satisfy the following functional
equation.

\begin{theorem}[Dyck path generating function, cf.~\cite{NT03b}]\label{theo:dyck}
The generating functions of Dyck paths $D({\bs u})$ and of arches 
$A({\bs u})$ satisfy
\begin{equation}\label{eqn:dyckfeq}
\begin{split}
D(u_0,\ldots,u_M)&=1+D(u_0,\ldots,u_M)A(u_0,\ldots,u_M),\\
A(u_0,u_1,\ldots,u_M)&=u_0 u_1\cdot\ldots\cdot u_M D(v_0,\ldots,v_M),
\end{split}
\end{equation}
where the monomials $v_k=v_k({\bs u})$ are given by
\begin{equation*}
v_0({\bs u})=u_0u_1^2u_2^2\cdot\ldots\cdot u_M^2, \qquad v_k({\bs u})
=\prod_{l=k}^M u_l^{\binom{l}{k}} \qquad (k=1,\ldots,M).
\end{equation*}
\qed
\end{theorem}

\noindent \textbf{Remarks.} \textit{i)} For $M=2$, the above theorem appears in 
\cite[Sec.~3.1]{NT03b}. For arbitrary $M$, it is (incorrectly) stated in \cite[p.~713]{NT03b}.
The proof for $M=2$ given in \cite{NT03b} generalises to arbitrary $M$. It uses
a last passage decomposition argument and the additivity of the counting parameters with
respect to the sequence construction, yielding the first equation in the theorem. 
Note that the arch decomposition of Dyck paths, as depicted in Figure \ref{fig:dyck},
implies the equivalent equation $D(\bs u)=1/(1-A(\bs u))$. 
For $n\in\mathbb N$, there is a bijection between the set of arches of length $2n$ and 
the set of Dyck paths of length $2n-2$, by identifying an arch, where the bottom layer 
has been removed, with the corresponding Dyck path. This implies, after a short calculation 
(see also \cite[Thm.~1]{R06}), the second equation in Theorem \ref{theo:dyck}.\\
\textit{ii)} The induced functional equation for the generating function $E(\bs u)=
D(\bs u)-1$ is an example of a $q$-functional equation as in Definition \ref{defi:qfunc},
with a square-root singularity in the generating function in the ``undeformed'' case 
$(u_1,\ldots,u_M)=(1,\ldots,1)$. See the following section and Section~\ref{sec:Dyckrev}.

\medskip

For a probabilistic interpretation of the counting parameters, as in the previous subsection,
let $p_{n_0,n_1,\ldots,n_M}$ denote the number of Dyck 
paths of length $2n_0$, where the $k$-th moment of height has the value $n_k$, 
for $k=1,\ldots,M$. We clearly have $0<\sum_{n_1,\ldots, n_M}p_{n_0,n_1,\ldots,n_M}<\infty$
for $n_0>0$. In the uniform ensemble, the random variables $\widetilde X_{k,n_0}$, which assign 
the $k$-th moment of height to a random Dyck path of length $2n_0$, have the joint distribution

\begin{equation}\label{eq:pxt}
\mathbb P(\widetilde X_{1,n_0}=n_1, \ldots, \widetilde X_{M,n_0}=n_M)=
\frac{p_{n_0,n_1,\ldots,n_M}}{\sum_{n_1,\ldots, n_M}p_{n_0,n_1,\ldots,n_M}}.
\end{equation}
The moments $\widetilde m_{k_1,\ldots,k_M}(n_0)$ of the joint distribution are given by

\begin{equation*}
\widetilde m_{k_1,\ldots,k_M}(n_0)=\frac{\sum_{n_1,\ldots, n_M} n_1^{k_1}
\cdot\ldots\cdot n_M^{k_M}p_{n_0,n_1,\ldots,n_M}}{\sum_{n_1,\ldots, n_M}
p_{n_0,n_1,\ldots,n_M}}.
\end{equation*}
In the previous subsection, we introduced the normalised
random variable eqn.~\eqref{form:DyckRV},
\begin{equation}\label{form:DyckRV2}
(X_{1,n_0},X_{2,n_0},\ldots, X_{M,n_0})= \left(\frac{\widetilde X_{1,n_0}}{n_0^{(1+2)/2}},
\frac{\widetilde X_{2,n_0}}{n_0^{(2+2)/2}},\ldots, 
\frac{\widetilde X_{M,n_0}}{n_0^{(M+2)/2}}\right).
\end{equation}
The normalised random variable $(X_{1,n_0},\ldots,X_{M,n_0})$ has moments 
$m_{k_1,\ldots,k_M}(n_0)$, given by
\begin{equation*}
m_{k_1,\ldots,k_M}(n_0)=\frac{\widetilde m_{k_1,\ldots,k_M}(n_0)}{
n_0^{(1+2)k_1/2+(2+2)k_2/2+\ldots+(M+2)k_M/2}}.
\end{equation*}
We argued above that these numbers should tend to a finite limit, as 
$n_0$ approaches infinity. As we will see below, this is indeed the 
case. Hence we may define
\begin{equation}\label{eqn:momlim}
m_{k_1,\ldots,k_M}=\lim_{n_0\to\infty}m_{k_1,\ldots,k_M}(n_0).
\end{equation}

A careful analysis of the functional equation eqn.~\eqref{eqn:dyckfeq}, which 
will be performed in the general case from Section \ref{sec:ana} onwards, 
yields the following result for the numbers $m_{k_1,\ldots,k_M}$. Its proof 
is deferred until Section~\ref{sec:Dyckrev}.

\begin{prop}\label{prop:dyckmom}
The normalised moments $m_{k_1,\ldots,k_M}$ of Dyck paths eqn.~\eqref{eqn:momlim}
are given by
\begin{equation}\label{eqn:dyckres}
\frac{m_{k_1,\ldots,k_M}}{k_1!\cdot\ldots\cdot k_M!}=
\frac{1}{f_{\bs0} u_c^{\gamma_{\bs k}-\gamma_{\bs 0}}}
\frac{\Gamma(\gamma_{\bs 0})}{\Gamma(\gamma_{\bs k})}
f_{\bs k},
\end{equation}
where $u_c=1/4$, where $\Gamma(z)$ denotes the Gamma function, and where the numbers 
$\gamma_{\bs k}$ and $f_{\bs k}$ are defined in eqn.~\eqref{form:exc2}. 
The moments have an entire exponential generating function. Hence, they define 
a unique random variable with moments $m_{\bs k}$.
\end{prop}

\noindent \textbf{Remarks.}
\textit{i)} Proposition~\ref{prop:dyckmom} implies convergence in distribution
and moment convergence of the normalised random variables eqn.~\eqref{form:DyckRV2}. 
This yields an alternative
proof of the convergence statement of Proposition~\ref{prop:reldyckexc}. However,
it does not establish a connection between the limit random variable and
the excursion moments, as in Proposition~\ref{prop:reldyckexc}.\\
\textit{ii)} Since the moments $m_{\bs k}$ define a unique random variable, the 
result for the excursion moments eqn.~\eqref{eqn:exmomexpl} follows from 
Proposition~\ref{prop:dyckmom} by Proposition \ref{prop:reldyckexc}. \\

\medskip

For the solution of a general $q$-functional equation (Definition \ref{defi:qfunc}), 
we obtain a similar result. Write $G(\bs u)=\sum_{n_1,\ldots,n_M}p_{n_0,n_1,\ldots,n_M}
u_0^{n_0}u_1^{n_1}\cdot\ldots\cdot u_M^{n_M}$ for such a solution. Under mild assumptions, 
the corresponding random variables in eqn.~\eqref{eq:pxt} are well-defined, and the same 
asymptotic analysis as above can be performed. We have the following theorem which is, 
in conjunction with Theorem \ref{theo:excmom}, the main result of our paper. Its proof is 
deferred until Section~\ref{sec:mom} and Section~\ref{sec:Dyckrev}.
\begin{theorem}[Limit distribution for $q$-functional equations]\label{theo:probdist}
Let a $q$-functional equation Definition \ref{defi:qfunc} with solution 
\begin{equation*}
G(u_0,u_1,\ldots,u_M)=\sum_{n_1,\ldots,n_M}p_{n_0,n_1,\ldots,n_M}u_0^{n_0}u_1^{n_1}
\cdot\ldots\cdot u_M^{n_M}
\end{equation*}
as in Proposition~\ref{prop:uni} be given, and let Assumption~\ref{ass} be satisfied.
Assume that the numbers $A_i$ in Proposition \ref{prop:rek} are
positive, i.e., $A_i>0$ for $i=0,\ldots,M-1$. Then, the random variables 
$(\widetilde X_{1,n_0},\ldots,\widetilde X_{M,n_0})$ eqn.~\eqref{eq:pxt} are
well-defined for almost all $n_0$. The following conclusions hold.

\begin{itemize}
\item[i)]
The numbers $m_{k_1,\ldots,k_M}$ eqn.~\eqref{eqn:momlim}, which are derived from the
moments of the normalised random variables $(X_{1,n_0},\ldots,X_{M,n_0})$ 
eqn.~\eqref{form:DyckRV2}, define a unique random variable
$(Y_1,\ldots,Y_M)$ with moments $m_{k_1,\ldots,k_M}$. We have convergence in distribution,
\begin{equation*}
(X_{1,n_0},\ldots, X_{M,n_0}) \stackrel{d}{\to} (Y_1,\ldots,Y_M)\qquad (n_0\to\infty),
\end{equation*}
and we have moment convergence.

\item[ii)]
The limiting random variable $(Y_1,\ldots,Y_M)$ is explicitly given by
\begin{equation*}
(Y_1,\ldots,Y_k)=(c_1X_1,\ldots,c_MX_M),
\end{equation*}
where the constants $c_k>0$ and the random variables $X_k$ are
\begin{equation*}
c_k =2^{\frac{k+2}{2}}\frac{2\mu_0\cdot\ldots\cdot \mu_{k-1}}{4^k k!}, \qquad
X_k=\int_{0}^1 e^k(t)\,{\rm d}t \qquad (k=1,\ldots,M).
\end{equation*} 
Here, the numbers $\mu_i>0$ are defined in Proposition \ref{prop:rek}, and $e(t)$ 
denotes a standard Brownian excursion of duration 1.

\end{itemize}
\end{theorem}

\noindent \textbf{Remarks.} \textit{i)} A recursion for the moments of the limit distribution
appears in Theorem~\ref{theo:excmom} above, see also Proposition \ref{prop:rek} below.\\
\textit{ii)} For Dyck paths counted by length and $k$-th moments of height, the 
above statements follow already from Proposition~\ref{prop:reldyckexc}. More 
generally, as was argued in the remark following 
Proposition~\ref{prop:reldyckexc}, the above result can be shown to hold for models 
of simply generated trees \cite{MM78}, counted by number of vertices and $k$-th moments 
of internal path length. This follows from the polynomial convergence of the 
depth first process derived from simply generated trees towards the Brownian excursion 
\cite{Bel72, A93, MM03, D03}. For a general $q$-functional equation, such a connection 
is not known to exist.\\
\textit{iii)} Our method of proof also allows to study corrections to the asymptotic behaviour, 
compare the discussion in Section~\ref{sec:dom}, and \cite{R02} for the case $M=1$. These 
cannot be obtained by the methods of \cite{MM03}. 

\subsection{Structure of the paper}

The remainder of the paper is organised as follows. In Section~\ref{sec:qfunc}, we 
introduce $q$-shifts (Definition~\ref{def:qshift}) and $q$-functional equations 
(Definition~\ref{defi:qfunc}). Our results rely on an application of the 
multivariate moment method (see e.g.~\cite[Ch.~6.1.]{JLR00} or \cite[Sec.~30]{Bill95}). 
Hence in Section~\ref{sec:gk}, eqn.~\eqref{form:facarmom}, we introduce 
\emph{factorial moment generating functions}, which are derivatives of the solution 
of the functional equation, evaluated at $u_1=\ldots=u_M=1$. We study their properties 
by a combinatorial analysis of derivatives of the functional equation, using a 
multivariate generalisation of Faa di Bruno's formula \cite{CS96}. In Section 
\ref{sec:ana}, we study the singular behaviour of the factorial moment generating 
functions, in the case of a square root singularity as the dominant singularity 
of the size generating function. Proposition~\ref{prop:rek} gives a recursion for 
the amplitudes, which describe the leading singular behaviour of the factorial moment 
generating functions.  Our method is also called \emph{moment pumping} \cite{FPV98}. 
We employ in Section~\ref{sec:dom} an alternative 
(rigorous) method to obtain the recursion, which originates from the 
\emph{method of dominant balance} of statistical mechanics \cite{PB95a, R02}. 
This method is generally easier to apply than an analysis of the functional equation 
via Faa di Bruno's formula, and yields an algorithm for obtaining the amplitude recursion.
In addition, the method allows to analyse corrections to the limiting behaviour. 
The behaviour of the moments follows then by standard methods 
from singularity analysis of generating functions \cite{FO90,FS07}.
In Proposition~\ref{theo:pde}, we give a quasi-linear partial differential equation 
for the generating function of the amplitudes. Growth estimates for the amplitudes 
(and hence for the moments) are obtained in Section~\ref{sec:growth}, by an analysis 
of the (singular) partial differential equation. Existence and uniqueness of a limit distribution 
is then guaranted by L\'evy's continuity theorem, see Section~\ref{sec:mom}. The connection 
to Brownian excursions follows in Section~\ref{sec:Dyckrev}, by a comparison with 
Dyck paths. Possible applications of our method are discussed in a concluding section.

\section{$q$-functional equations}\label{sec:qfunc}

Let $\mathbb C [[\bs u]]$ denote the ring of formal power series with complex coefficients 
in the (commutative) variables ${\bs u}=(u_0,u_1,\ldots, u_M)$. Let $\mathbb C [\bs u]$ 
denote the ring of polynomials with complex coefficients in the variables ${\bs u}$. We 
set ${\bs u}_0=(u_0,1,\ldots,1)$. For ${\bs u}=(u_0,u_1,\ldots,u_M)$ and ${\bs n}\in
\mathbb N_0^{1+M}$, we define $\bs u^{\bs n}$ to be the monomial ${{\bs u}^{\bs n}}=
u_0^{n_0}\cdot\ldots\cdot u_M^{n_M}$. We employ the notation ${\bs u}_+=(u_1,\ldots,u_M)$, 
the plus sign indicating that the first component of $\bs u$ is omitted. For $n\in\mathbb 
N$ and $r>0$, let $D^n_r$ denote the open polydisc 
\begin{equation*}
D^n_r=\{ (x_1,\ldots, x_n)\in\mathbb C^n: |x_k|<r \;\mbox{for all}\; k=1,\ldots, n\}.
\end{equation*}
Let $\mathcal{H}_r(\bs x)$ denote the ring of power series in $\bs x=(x_1,\ldots,x_n)$ 
with complex coefficients, which are convergent in $D^n_r$. 

\begin{defi}[$q$-shift]\label{def:qshift}
Fix $M\in\mathbb N$. Let for $k=0,\ldots,M$ formal power series 
$v_k(\bs u)\in\mathbb C[[\bs u]]$ be given, and write
$\bs v=\bs v(\bs u)=(v_0(\bs u),
v_1(\bs u),\ldots,v_M(\bs u))$. Assume that there is a
number $d=d(\bs v)$ satisfying $1<d\le\infty$, such that
\begin{equation*}
v_0(\bs u)\in \mathcal{H}_d(u_1,\ldots,u_M)[[u_0]],\qquad
v_k(\bs u)\in \mathcal{H}_d(u_k,\ldots,u_M) \qquad (k=1,\ldots, M).
\end{equation*}
Assume that $\bs v(\bs u)$ satisfies
\begin{equation*}
{\bs v({\bs u}_0)={\bs u}_0}, \qquad
v_0(0,{\bs u}_+)\equiv0, \qquad \frac{\partial v_k}
{\partial u_k}({\bs u}_0)\equiv1 \qquad (k=0,\ldots, M).
\end{equation*} 
Let $r=r(\bs v)$ be a number $0<r\le\infty$ such that
\begin{equation*}
\left\{(v_1(\bs u_+),\ldots, v_M(\bs u_+) ): \bs u_+\in D_d^M\right\}\subseteq D_r^M.
\end{equation*}
Then, $\bs v$ is called a \emph{$q$-shift}, $d(\bs v)$ is called
the \emph{domain of $\bs v$}, and $r(\bs v)$ is called the 
\emph{range of $\bs v$}.

\end{defi}

\noindent \textbf{Remarks.}
\textit{i)} For $M=1$, a simple example of $q$-shift appears in 
the $q$-difference equation eqn.~\eqref{form:geneq2}, with
$v_0(u_0,u_1)=u_0u_1$ and $v_1(u_1)=u_1$. This motivates 
the name for the generalisation. \\
\textit{ii)} For a $q$-shift ${\bs v}$, its $k$-th component 
$v_k(\bs u)$ does not depend on $u_l$, where $l=0,\ldots,k-1$,
and it does depend on $u_k$. Since ${\bs v({\bs u}_0)={\bs u}_0}$,
one may interpret $\bs v(\bs u)$ for $\bs u\ne\bs u_0$
as a ``deformation'' of $\bs u_0$. The condition $v_0(0,{\bs u}_+)\equiv0$
is imposed to ensure that composition of formal power series is well-defined, 
see below.\\
\textit{iii)}
The identity $id: {\bs u}\mapsto {\bs u}$ is a $q$-shift with infinite domain and 
range. For two $q$-shifts ${\bs v}$ and ${\bs w}$, their composition 
${\bs v}\circ {\bs w}$ is well-defined if $r(\bs w)\le d(\bs v)$. As is readily checked, 
${\bs v}\circ {\bs w}$ is a $q$-shift in that case, with domain $d({\bs v}\circ {\bs w})=d(\bs w)$ 
and range $r({\bs v}\circ {\bs w})=r(\bs v)$. If for a $q$-shift ${\bs v}$  we have 
$r(\bs v)\le d(\bs v)$, we can thus consider \emph{iterated} $q$-shifts ${\bs v}^{[n]}$, where
\begin{equation*}
{\bs v}^{[0]}=id, \qquad {\bs v}^{[n]}={\bs v}\circ{\bs v}^{[n-1]} \qquad (n\in\mathbb N).
\end{equation*}
The power series ${\bs v}^{[n]}$ is a $q$-shift for $n\in \mathbb N_0$.\\
\textit{iv)} An important subclass (with infinite domain and range) are \emph{polynomial} 
$q$-shifts, i.e., $q$-shifts $\bs v$ satisfying $v_k(\bs u) \in 
\mathbb C[\bs u]$ for $k=0,\ldots, M$. Examples are given by 
the monomial $q$-shifts
\begin{equation*}
v_k({\bs u}) = {\bs u}^{{\bs n}_k},
\end{equation*}
where ${\bs n}_k\in\mathbb N_0^{1+M}$, $({\bs n}_k)_k=1$ 
and $({\bs n}_k)_l=0$ for $l=0,\ldots,k-1$ and $k=0,\ldots,M$. These
appear in eqns.~\eqref{form:geneq2} and \eqref{form:geneq3}, and for 
Dyck paths in eqn.~\eqref{eqn:dyckfeq}. In these examples, 
we have $({\bs n}_l)_{l+1}\ne 0$ for $l=0,\ldots, M-1$.

\medskip

For a formal power series $G({\bs u})\in \mathcal{H}_d(\bs u_+)[[u_0]]
\subseteq\mathbb C [[\bs u]]$ and a $q$-shift ${\bs v}$ satisfying
$r(\bs v)\le d$, we define $H({\bs u}) \in \mathcal{H}_{d(\bs v)}(\bs u_+)[[u_0]]$ by
\begin{equation*}
H({\bs u})= G({\bs v}({\bs u})).
\end{equation*}
This is well-defined since $\mathcal{H}_d({\bs v}_+({\bs u}))\subseteq 
\mathcal{H}_{d(\bs v)}(\bs u_+)$ and $\mathbb C[[v_0({\bs u})]]\subseteq
\mathcal{H}_{d(\bs v)}(\bs u_+)[[u_0]]$, due to $v_0(0,{\bs u}_+)\equiv0$. 
We are interested in derivatives of $H({\bs u)}$. For 
clarity of presentation, we will use the multi-index notation, and 
Greek indices $\bs\mu, \bs\nu,\bs\rho$ will denote vectors with 
non-negative integer entries. For $F({\bs u}) \in\mathbb C [[\bs u]]$ and 
$\bs\nu=(\nu_0,\ldots,\nu_M)\in
\mathbb N_0^{1+M}$, we write the derivative $F_{\bs \nu}({\bs u}) 
\in\mathbb C [[\bs u]]$ of $F(\bs u)$ of order $\bs\nu$ as
\begin{equation*}
F_{\bs\nu}({\bs u}) = \partial_0^{\nu_0}\cdots 
\partial_M^{\nu_M}F({\bs u}),\qquad \partial_k^{\nu_k}=
\frac{\partial^{\nu_k}}{\partial u_k^{\nu_k}}\qquad (k=0,\ldots,M),
\end{equation*}
where we use the convention that $F_{\bs 0}(\bs u)=F(\bs u)$.
If $G({\bs u}) \in \mathcal{H}_d(\bs u_+)[[u_0]]\subset\mathbb 
C[[\bs u]]$, we also have $G_{\bs\nu}({\bs u})\in\mathcal{H}_d(\bs u_+)
[[u_0]]$. Set $|\bs\nu|=\nu_0+\ldots+\nu_M$ and $\bs\nu!=\nu_0!
\cdot\ldots\cdot \nu_M!$. We define a total order in $\mathbb N_0^{1+M}$
as follows.
 
\begin{defi}[total order $\prec$]\label{def:totord}
For $\bs\mu =(\mu_0,\ldots,\mu_M)\in\mathbb N_0^{1+M}$ and 
$\bs\nu=(\nu_0,\ldots,\nu_M)\in\mathbb N_0^{1+M}$, we write 
$\bs\mu \prec\bs\nu$ if either $|\bs\mu |<|\bs\nu|$, or if $|\bs\mu |=|\bs\nu|$, 
there exists an index  $k\in\{0,\ldots,M\}$, such that
$\mu_k>\nu_k$ and $\mu_i=\nu_i$ for $i=0,\ldots,k-1$.
\end{defi}

The total order introduced above will be used to label terms appearing
in derivatives of $H({\bs u)}$. We have the following lemma.  For 
$i\in\{0,\ldots,M\}$, let ${\bs e}_i\in\mathbb C^{1+M}$ denote the 
unit vector in direction $i$, given by $({\bs e}_i)_k=\delta_{i,k}$ for $k=0,\ldots,M$.

\begin{lemma}\label{lemma:H}
Let $G({\bs u}) \in \mathcal{H}_d(\bs u_+)[[u_0]]$, and let a 
$q$-shift $\bs v$ with domain $d(\bs v)$ and range $r(\bs v)\le d$ be given. 
Set $H({\bs u)}=G({\bs v(\bs u)})\in \mathcal{H}_{d(\bs v)}(\bs u_+)[[u_0]]$ 
and fix $\bs\nu\ne {\bf 0}$.
\begin{itemize}

\item[i)] 
We have $H_{\bs\nu}({\bs u})\in\mathcal{H}_{d(\bs v)}(\bs u_+)[[u_0]]$.
For every $\bs \mu$ satisfying $\mathbf{0}\ne\bs\mu \preceq\bs\nu$, there
exists a coefficient $A(\bs\mu,\bs u)\in\mathcal{H}_{d(\bs v)}(\bs u_+)[[u_0]]$, 
independent of the choice of $G({\bs u})$, such that
\begin{equation}\label{form:chain}
H_{\bs\nu}({\bs u}) = 
\sum_{{\bf 0}\ne\bs\mu \preceq\bs\nu}
G_{\bs\mu }({\bs v}({\bs u})) \cdot A(\bs\mu ,{\bs u}).
\end{equation}

\item[ii)] The coefficient $A(\bs\mu,{\bs u})$ in eqn.~\eqref{form:chain}
is, for $\bs\mu=\bs\nu$, given by
\begin{equation}\label{form:anuu}
A(\bs\nu,{\bs u})=\prod_{k=0}^M \left( \frac{\partial v_k}{\partial u_k}
({\bs u}) \right)^{\nu_k}.
\end{equation}
In particular, we have $A(\bs\nu,{\bs u_0})=1$. 

\item[iii)]
Fix $i\in\{0,\ldots,M-1\}$. If $\nu_{i+1}>0$, 
we have $\bs\nu-{\bs e}_{i+1}+{\bs e}_i \prec \bs\nu$.
The coefficient $A(\bs\mu,{\bs u})$ in eqn.~\eqref{form:chain}
is, for $\bs\mu=\bs\nu-{\bs e}_{i+1}+{\bs e}_i$ and 
$\bs u=\bs u_0$, given by
\begin{equation*}
A(\bs\nu-{\bs e}_{i+1}+{\bs e}_i,{\bs u_0})= 
\nu_{i+1}\frac{\partial v_i}{\partial u_{i+1}}({\bs u_0}).
\end{equation*}

\item[iv)] 
For real numbers $r,r_0,\ldots,r_M$, where $r_{k+1}>r_k$ for 
$k=0,\ldots,M-1$, define $\alpha_{\bs\mu }=r+\sum_{k=0}^M r_k\mu_k$. 
Then, for indices $\bs\mu \prec\bs\nu$ in eqn.~\eqref{form:chain} 
satisfying $A(\bs\mu ,{\bs u})\not\equiv0$, we have
$\alpha_{\bs\mu }<\alpha_{\bs\nu}$.
\end{itemize}
\end{lemma}

\noindent {\bf Remarks.} \textit{i)} The coefficient $A(\bs\mu ,{\bs u})$ 
in eqn.~\eqref{form:chain} might be chosen to vanish. E.g., fix $\bs\nu=(1,1)$ 
and consider $\bs\mu =(0,1)$. We have $\bs\mu\prec\bs
\nu$ but might choose $A(\bs\mu ,{\bs u})\equiv0$, as is readily verified by an 
explicit calculation, using the $q$-shift property $v_k({\bs u}) \in 
\mathbb C[[u_k,\ldots,u_M]]$. \\
\textit{ii)} The property \textit{iv)} will be used for exponent estimates
in Proposition~\ref{prop:est} below.

\begin{proof}
{\it i)} This is an application of the chain rule. Successive differentiation leads to
\begin{equation}\label{form:succ}
H_{\bs\nu}({\bs u}) = 
\sum_{0\le|\bs\mu |\le|\bs\nu|}
G_{\bs\mu }({\bs v}({\bs u})) \cdot A(\bs\mu ,{\bs u}),
\end{equation}
where $A(\bs\mu ,{\bs u})$ is independent of $G({\bs u})$. To analyse the 
effect of the $q$-shift property $v_k({\bs u})\in\mathbb {\cal H}_{d(\bs v)}(u_k,\ldots, u_M)$, 
consider for $k\in\{0,\ldots, M\}$ the equation
\begin{equation}\label{form:d1}
\frac{\partial H}{\partial u_k}({\bs u}) =
\sum_{l=0}^M \frac{\partial G}{\partial v_l}({\bs v}({\bs u}))\frac{\partial v_l}
{\partial u_k}({\bs u})=\sum_{l=0}^k \frac{\partial G}{\partial v_l}({\bs v}
({\bs u}))\frac{\partial v_l}{\partial u_k}({\bs u}).
\end{equation}
It shows that derivatives of $H({\bs u})$ w.r.t.~$u_k$ do not contribute to derivatives of 
$G({\bs v}({\bs u}))$ w.r.t.~$v_{k+1},\ldots,v_M$. A similar statement holds for 
higher derivatives. This implies that the numbers $\nu_k$ in eqn.~\eqref{form:succ}
can contribute to the numbers $\mu_0,\ldots,\mu_k$ only. If $|\bs\mu |=|\bs\nu|$, we thus 
have  $\nu_k\le \mu_0+\ldots+\mu_k$ for $k\in\{0,\ldots,M\}$. We show that $\bs\mu 
\preceq\bs\nu$. Assume w.l.o.g. that $|\bs\mu |=|\bs\nu|$. If $\mu_0>\nu_0$, we have  
$\bs\mu \prec \bs\nu$, and the claim follows. Otherwise, we have $\mu_0=\nu_0$. 
Thus, $\nu_k$ does not contribute to $\mu_0$ for $k\in\{1,\ldots,M\}$. The previous 
argument yields that $\nu_k\le \mu_1+\ldots+\mu_k$ for $k\in\{1,\ldots,M\}$. Now repeat 
the above argument until $\mu_0=\ldots=\mu_{M-1}=\nu_0=\ldots=\nu_{M-1}$. Then 
$\mu_M=\nu_M$ since $|\bs\mu |=|\bs\nu|$. Thus $\bs\mu =\bs\nu$, and the statement 
is shown. 

\noindent {\it ii)}
This is seen by an explicit calculation using eqn.~\eqref{form:d1}. Successively applying 
derivatives w.r.t. $u_0$, $u_1$, $\ldots,u_M$ yields eqn.~\eqref{form:anuu}. Together 
with the $q$-shift property ${\partial v_k}/{\partial u_k}({\bs u_0})=1$, it follows that 
$A(\bs\nu,{\bs u_0})=1$. Statement \textit{iii)} is shown by an analogous calculation.

\noindent {\it iv)} Consider first the case $|\bs\mu |=|\bs\nu|$. Note that, by differentiation, 
the numbers $\nu_k$ do not contribute to $\mu_{k+1},\ldots,\mu_M$, for $k\in\{0,\ldots,M\}$. 
Thus, the contribution of $\nu_k$ to $\alpha_{\bs\mu }$ is maximal if $\mu_k=\nu_k$ for 
$k=0,\ldots,M$.  We get $\bs\mu =\bs\nu$ and $\alpha_{\bs\mu }=\alpha_{\bs\nu}$. Since 
this maximum is unique, $\bs\mu \prec \bs\nu$ implies $\alpha_{\bs\mu }<\alpha_{\bs\nu}$.
Let now $|\bs\mu |<|\bs\nu|$. For $k\in\{0,\ldots,M\}$, denote by ${\widetilde \nu}_k$ the 
number of derivatives w.r.t.~$u_k$, which contribute to $\bs\mu $. Set $\widetilde{\bs\nu}=
({\widetilde \nu}_0, \ldots, {\widetilde \nu}_M)$. Then clearly $|\bs\mu |=|\widetilde{\bs\nu}|$ 
and $\alpha_{\widetilde{\bs\nu}}<\alpha_{\bs\nu}$. The reasoning for the case $|\bs\mu |=
|\bs\nu|$ can now be applied to the present case, with $\bs\nu$ replaced by $\widetilde{\bs\nu}$.
This yields $\alpha_{\bs\mu }\le\alpha_{\widetilde{\bs\nu}}$. Thus $\alpha_{\bs\mu }<
\alpha_{\bs\nu}$, and the statement is shown.
\end{proof}

\begin{defi}[$q$-functional equation]\label{defi:qfunc}
Let $P(\bs u, y_1,\ldots, y_N)$ be a formal power series $P(\bs u,
\bs y)\in\mathcal{H}_{d}(\bs u_+)[[u_0,\bs y]]$, for
a number $d$ satisfying $1<d\le\infty$. Assume that 
$P(0, {\bs u}_+, \bs 0)\equiv 0$, and that $\frac{\partial P}{\partial y_j}(0, 
{\bs u}_+,\bs 0) \equiv0$ for $j=1,\ldots,N$. Let ${\bs v^{(\it j)}}$ be a $q$-shift
with domain $d(\bs v^{(j)})$ and range $r(\bs v^{(j)})\le d$, for $j=1,\ldots,N$.

\begin{itemize}

\item[i)]

If the above assumptions are satisfied, the equation
\begin{equation}\label{form:funceq}
G({\bs u})=P(\bs u, H^{(1)}({\bs u)},\ldots, H^{(N)}
({\bs u)}),
\end{equation}
where $H^{(j)}({\bs u)}=G({\bs v^{(\it j)}(\bs u)})$ for $j=1,\ldots, N$,
is called a {\em $q$-functional equation}.
\item[ii)]
Let the above assumptions be satisfied. If there is a number $d(G)$ such that
\begin{equation*}
1< d(G)\le \min_{1\le j\le N}\{d(\bs v^{(j)})\},
\end{equation*}
and a formal power series $G(\bs u)\in\mathcal{H}_{d(G)}(\bs 
u_+)[[u_0]]$ satisfying eqn.~\eqref{form:funceq}, then $G(\bs u)$ 
is called a \emph{solution} of the $q$-functional equation.
\end{itemize}
\end{defi}

\noindent \textbf{Remarks.} 
\textit{i)} An example of a $q$-functional equation for $M=1$ is 
given by the $q$-difference equation eqn.~\eqref{form:geneq2}. 
This motivates the name for the generalisation. Examples for $M>1$ 
appear in eqn.~\eqref{form:geneq3}. Examples of $q$-functional equations 
frequently satisfy $P(\bs u, \bs y)\in\mathbb C[\bs u, \bs y]$ 
with $d=\infty$, and with monomial $q$-shifts. We infer from the 
functional equation eqn.~\eqref{eqn:dyckfeq} for Dyck paths that 
the power series $E({\bs u})=D({\bs u})-1$ satisfies the 
$q$-functional equation
\begin{equation}\label{form:exfunc}
E({\bs u}) = u_0u_1\cdot\ldots\cdot u_M\left( E({\bs u})+1\right)
\left( E({\bs v(\bs u)})+1\right).
\end{equation}
Specialising to ${\bs u}={\bs u_0}$ yields a quadratic equation, which
can be explicitly solved for $E({\bs u_0})$. We get the well-known
result
\begin{equation*}
E({\bs u_0})=\frac{1-2u_0-\sqrt{1-4u_0}}{2u_0}.
\end{equation*}

\noindent \textit{ii)} If $Q(u,y)\in\mathbb 
C[[u,y]]$ satisfies $Q(0,0)=0$ and $\frac{\partial Q}{\partial y}(0,y)
\not\equiv1$, a formal power series $G(u)\in\mathbb C[[u]]$ such that 
$G(0)=0$ is uniquely defined as the solution of the equation 
$G(u)=Q(u,G(u))$. When studying such 
equations, we may assume without loss of generality that $\frac{\partial Q}
{\partial y}(0,y)\equiv0$. The above definition reduces to this setup, when 
restricting to $\bs u=\bs u_0$. \\

\medskip

A result about solutions of $q$-functional equations is given by the following proposition.
We use the vector notation ${\bs H(\bs u)}=(H^{(1)}({\bs 
u)},\ldots, H^{(N)}({\bs u)})$. 

\begin{prop}\label{prop:uni}
The $q$-functional equation of Definition~\ref{defi:qfunc},
\begin{equation}\label{form:func}
G({\bs u})=P(\bs u, \bs H(\bs u)),
\end{equation}
has a unique solution $G({\bs u}) \in \mathcal{H}_{d(P)}(\bs u_+)[[u_0]]$
satisfying $G(0,\bs u_+)\equiv0$, where $d(P)=\min_{1\le j\le N}\{d(\bs v^{(j)})\}$.
\end{prop}

\begin{proof}
If $G({\bs u})=\sum_{n_0\ge1} p_{n_0}({\bs u}_+)u_0^{n_0}$, then
$G({\bs v}^{(\it j)}({\bs u}))=\sum_{n_0\ge1} p_{n_0}({\bs v}_+^{(\it j)}
({\bs u}))\left(v_0^{(\it j)}({\bs u})\right)^{n_0}$. For $n_0\ge1$ fixed, we take the coefficient 
of $u_0^{n_0}$ in eqn.~\eqref{form:func}. Due to $v_0^{(j)}(0,{\bs u}_+)
\equiv0$ for $j=1,\ldots,N$ and the assumptions on the derivatives of $P({\bs y},{\bs u})$, we 
get the expression
\begin{equation}\label{form:Wr}
p_{n_0}({\bs u}_+)=W_{n_0}\left({\bs u}_+, \left\{ p_1({\bs v}_+^{(\it j)}
({\bs u})), \ldots,p_{n_0-1}({\bs v}_+^{(\it j)}({\bs u}) )\right\}_{j=1}^N
\right),
\end{equation}
for a power series $W_{n_0}\left({\bs u}_+, \bs p\right)\in\mathcal{H}_{d(P)}
(\bs u_+)\left[\bs p\right]$ in $M+N(n_0-1)$ variables. Thus $p_{n_0}
({\bs u}_+)$ is determined recursively in terms of $p_l({\bs u}_+)$, 
where $l=1,\ldots, n_0-1$. The recursion also shows that $p_{n_0}({\bs 
u}_+)\in \mathcal{H}_{d(P)}(\bs u_+)$. Thus $G({\bs u}) \in \mathcal{H}_{d(P)}
(\bs u_+)[[u_0]]$, and the proposition is proved.
\end{proof}

\noindent \textbf{Remarks.} 
\noindent \textit{i)} For a solution $G(\bs u)$ of a $q$-functional 
equation satisfying $G(\bs0)=0$, we will always assume $d(G)=d(P)$ in the following.\\
\textit{ii)} For the solution of a simple $q$-functional 
equation, an explicit expression may be given. This is e.g. the
case for some $q$-difference equations eqn.~\eqref{form:geneq2}, with $P(x,y)$
linear in $y$, see e.g.~\cite{Bou96}, or with $P(x,y)$ quadratic in $y$, see \cite{PB95a}. \\
\textit{iii)} It follows from $G({\bs u}) \in \mathcal{H}_{d(P)}
(\bs u_+)[[u_0]]$ that $G({\bs u_0})\in\mathbb C[[u_0]]$. For derivatives of 
$G({\bs u})$, which are also elements of $\mathcal{H}_{d(P)}(\bs u_+)[[u_0]]$, the 
same conclusion holds. Thus, the derivatives
\begin{equation}\label{form:facarmom}
g_{\bs\nu}(u_0) :=\left.\frac{1}{\bs\nu!}
G_{\bs\nu}({\bs u})\right|_{\bs u=\bs u_0}
\end{equation}
are formal power series, i.e., $g_{\bs\nu}(u_0)\in\mathbb C[[u_0]]$. They 
are called \emph{factorial moment generating functions}, for reasons to 
be explained in section \ref{sec:mom}.

\section{Factorial moment generating functions}\label{sec:gk}

Due to the following proposition, the factorial moment generating functions 
can be computed recursively, by successively differentiating the $q$-functional 
equation.

\begin{prop}\label{prop:solv}
Let a $q$-functional equation eqn.~\eqref{form:funceq} with solution 
$G(\bs u)$ as in Proposition~\ref{prop:uni} be given. Consider the derivative of order
$\bs\nu\ne {\bf 0}$ of the $q$-functional equation, 
evaluated at $\bs u=\bs u_0$. It is linear in 
$g_{\bs\nu}(u_0)$. Its r.h.s.~is a polynomial in 
$\{g_{\bs\mu }(u_0):\mathbf{0}\neq\bs
\mu \preceq \bs\nu\}$, with coefficients in 
$\mathbb C[[u_0,g_{\bs 0}(u_0)]]$.
\end{prop}

In order to prove Proposition \ref{prop:solv}, we analyse partial 
derivatives of eqn.~\eqref{form:funceq}. To this end, we employ a generalization 
of Faa di Bruno's formula \cite{CS96}, adapted to our situation. The following 
lemma will also be used for exponent estimates in the next section. For 
$k\in\{1,\ldots,K\}$, let $f_k(\bs u)\in\mathbb C[[\bs u]]$, and 
define ${\bs f(\bs u)}=(f_1({\bs u}), \ldots,
f_K({\bs u}))$. For $\bs\mu \in\mathbb N_0^{1+M}$, we 
use the notation ${\bs f(\bs u)}_{\bs\mu} =
((f_1)_{\bs\mu}({\bs u}),\ldots,(f_K)_{\bs\mu}
({\bs u}))$ for the derivative of $\bs f(\bs u)$ of order $\bs \mu$.

\begin{lemma}{\rm (cf.~\cite[Thm.~2.1]{CS96})}\label{lemma:faa}
Let a $q$-functional equation eqn.~\eqref{form:funceq} with solution 
$G(\bs u)$ as in Proposition~\ref{prop:uni} be given.
Its derivative of order $\bs\nu\ne\bs0$ satisfies
\begin{equation}\label{form:faa}
G_{\bs\nu}({\bs u})=\sum_{1\le|\bs\lambda|\le |\bs\nu|}
P_{\bs\lambda}(\bs u, \bs H(\bs u))\sum_{s=1}^{|\bs\nu|}
\sum_{p_s(\bs\nu,\bs\lambda)}
(\bs\nu!) \prod_{j=1}^s\frac{[(\bs u, \bs H(\bs u))_{\bs\mu _j}]^{\bs\kappa_j}}
{(\bs\kappa_j!)[\bs\mu _j!]^{|\bs\kappa_j|}},
\end{equation}
where the vectors $\bs\lambda$ and $\bs\kappa_j$ are $1+M+N$-dimensional, 
the vectors $\bs\mu _j$ are $1+M$-dimensional, for $j\in\{1,\ldots,s\}$, 
and the summation ranges over
\begin{equation*}
\begin{split}
p_s(\bs\nu,\bs\lambda)=\{&(\bs\kappa_1,\ldots, \bs\kappa_s; 
\bs\mu _1,\ldots, \bs\mu _s): |\bs\kappa_i|>0,\\
&{\bf 0}\lhd \bs\mu _1\lhd \ldots \lhd \bs\mu _s, 
\sum_{i=1}^s \bs\kappa_i=\bs\lambda \mbox{ and } \sum_{i=1}^s |\bs\kappa_i| 
\bs\mu _i= \bs\nu \}.
\end{split}
\end{equation*}
In the above equation, the total order $\lhd$ is defined as 
$\bs\mu \lhd\bs\nu$ if either $|\bs\mu |<|\bs\nu|$, or
if $|\bs\mu |=|\bs\nu|$, there exists an index $k\in\{0,\ldots,M\}$ such that
$\mu_k<\nu_k$ and $\mu_i=\nu_i$ for $i=0,\ldots,k-1$. 
\qed
\end{lemma}

\noindent \textbf{Remark.} If $M=1$, we have $(\mu_0,\mu_1)\lhd(\nu_0,\nu_1)$
if and only if $(\mu_1,\mu_0)\prec(\nu_1,\nu_0)$. The analogous statement for
higher values of $M$ is not true.

\begin{proof}[Proof of Proposition \ref{prop:solv}.]
For given $\bs\nu\neq\bs0$, we analyse the values $\bs\mu _j$ appearing in
Lemma \ref{lemma:faa}. We show that $|\bs\mu _j|\ge|\bs\nu|$ for
some value $j$, where $1\le j\le s \le |\bs\nu|$, implies $s=1$ and 
$\bs\mu _1=\bs\nu$. Together with Lemma \ref{lemma:H}, this implies that 
$\bs\mu \preceq \bs\nu$. The linearity will follow from an explicit 
calculation of the term containing $\bs\mu _1=\bs\nu$.

\smallskip

Assume that $|\bs\mu _j|\ge|\bs\nu|$ for some $j$, where 
$1\le j\le s\le |\bs\nu|$. The explicit form of 
$p_s(\bs\nu,\bs\lambda)$ in Lemma \ref{lemma:faa}
states that $\bs\nu=\sum_{i=1}^s |\bs\kappa_i| \bs\mu _i$. This implies 
that $|\bs\nu|=\sum_{i=1}^s |\bs\kappa_i| |\bs\mu _i|$. According 
to the assumption, this leads to $s=1$,  $|\bs\kappa_1|=1$ and 
$|\bs\mu _1|=|\bs\nu|$. This, in turn, implies that $\bs\mu _1=\bs\nu$.

Clearly, the r.h.s.~of eqn.~\eqref{form:faa}, when specialised to $\bs u=\bs u_0$,
is a polynomial in $g_{\bs\mu }(u_0)$ for $\mathbf{0}\neq\bs\mu\preceq\bs\nu$,
with coefficients in $\mathbb C[[u_0,g_{\bs 0}(u_0)]]$. 
To show linearity in $G_{\bs\nu}({\bs u_0})$, we note that 
the possible values of $\bs\kappa_1$ are $\bs\kappa_1={\bs e}_j$, 
where $j\in\{1,\ldots,1+M+N\}$. The sum of the terms with $|\bs\mu _1|=|\bs\nu|$
in the r.h.s.~of eqn.~\eqref{form:faa} is given by
\begin{equation*}
\sum_{j=1}^{1+M+N} P_{{\bs e}_j}({\bs u, \bs H(\bs u)})
[({\bs u, \bs H(\bs u)})_{\bs\nu}]^{{\bs e}_j}.
\end{equation*}
We now extract terms containing $G_{\bs\nu}({\bs u})$ from this 
expression, using Lemma \ref{lemma:H}, and group them to the l.h.s.~of 
eqn.~\eqref{form:faa}. In the resulting equation, the l.h.s.~$L(\bs u)$, 
when specialised to $\bs u=\bs u_0$, is given by
\begin{equation}\label{form:pref}
L(\bs u_0) =\left(1- \sum_{j=1}^{N} \frac{\partial P}{\partial y_j}
({\bs u_0, \bs G(\bs u_0)})\right) G_{\bs\nu}({\bs u_0}).
\end{equation}
Due to the assumptions on $P(\bs u,\bs y)$, the prefactor
of $G_{\bs\nu}({\bs u_0})$ in the above equation is not
identically vanishing. Thus, $G_{\bs\nu}({\bs u_0})$ is contained linearly in 
eqn.~\eqref{form:faa}, specialised to $\bs u=\bs u_0$.
\end{proof}

\section{Analytic generating functions}\label{sec:ana}

Let $Q(u,y)=P(u,1,\ldots,1,y,\ldots,y)$. The power series $G({\bs u}_0)$
satisfies the equation
\begin{equation}\label{pgfeqn}
G({\bs u}_0)=Q(u_0,G({\bs u}_0)).
\end{equation}
In the following, we specialise the class of $q$-functional equations. We are 
interested in the case  where $G({\bs u}_0)$ is analytic at $u_0=0$, with 
a square root as dominant singularity. This situation is generic for 
combinatorial constructions, see \cite[Thm.~10.6]{O95}, \cite[Prop.~1]{D97}, or
\cite[Sec.~7.4]{FS07}. Throughout the remainder of the article, we employ 
the following assumption.

\begin{ass}\label{ass}
Let a $q$-functional equation \eqref{form:funceq} as in Definition~\ref{defi:qfunc} 
be given. Let numbers $r,s$ such that $0<r,s\le\infty$ be given. Denote by
$\mathcal{H}_{r,d,s}(\bs u,\bs y)$ the ring of power series in 
$(\bs u,\bs y)=(u_0,\bs u_+,\bs y)$, which
are convergent if  $|u_0|<r$, if $|u_k|<d$ for $k=1,\ldots, M$, and if $|y_k|<s$ for 
$k=1,\ldots,N$. In addition to the assumptions in Definition~\ref{defi:qfunc},
we assume the following properties.
\begin{itemize}
\item[{\it i)}] $P({\bs u},{\bs y})\in\mathcal{H}_{r,d,s}(\bs u, 
\bs y)$, and its Taylor coefficients are non-negative real numbers,
when expanded about $(\bs u,\bs y)=(\bs 0,\bs 0)$.
For $j=1,\ldots,N$, each coordinate series of the $q$-shift 
${\bs v}^{(j)}({\bs u})$ has non-negative
coefficients only, if expanded about $\bs u=\bs 0$.

\item[{\it ii)}] Let $Q(u,y)=P(u,1,\ldots,1,y,\ldots,y)$. There exist numbers $(u_c,y_c)$
satisfying $0<u_c<r$ and $0<y_c<s$, such that 
\begin{equation}\label{form:lett}
\begin{split}
Q(u_c,y_c)=&y_c, \qquad \left.\frac{\partial Q}{\partial y}(u_c,y)\right|_{y=y_c}=1, \\
B:=\frac{1}{2}\left.\frac{\partial ^2 Q}{\partial y^2}(u_c,y)\right|_{y=y_c}&>0, \qquad 
C:=\left.\frac{\partial Q}{\partial u}(u,y_c)\right|_{u=u_c}>0.
\end{split}
\end{equation}
\item[{\it iii)}] The solution $G(\bs u)$ of the $q$-functional equation as in 
Proposition~\ref{prop:uni} has the property that $G({\bs u_0})=\sum_{n\ge 1} 
p_n u_0^n$ is \emph{aperiodic}, i.e., there exist indices $1\le i<j<k$ such that 
$p_ip_jp_k\ne0$, while $\gcd(j-i,k-i)=1$.
\end{itemize}
\end{ass}

\noindent {\bf Remarks.} \textit{i)} When restricting to $\bs u=\bs u_0$,
the above assumptions (together with the assumptions in Definition~\ref{defi:qfunc}) 
reduce to the setup for implicitly defined power series which is usual in enumerative 
combinatorics, see \cite[Thm.~10.6]{O95}, \cite[Prop.~1]{D97}, and \cite[Sec.~7.4]{FS07}. \\
\textit{ii)} Assumption \ref{ass} {\it i)} implies that the coefficients $p_{\bs n}$ in 
$G({\bs u})= \sum_{\bs n\ge 0} p_{\bs n}{\bs u^{\bs n}}$ 
are non-negative. For combinatorial constructions, such positivity assumptions are common. 
However, systems of functional equations, arising from a combinatorial construction 
with positive coefficients, might be reduced to an equation of the above form 
having \emph{negative} coefficients. In that situation, other types of singularity might 
appear. This is, e.g., the case for discrete meanders \cite{NT03}. The methods used in 
this paper can be adapted to treat such cases. $P({\bs u},
{\bs y})\in\mathcal{H}_{r,d,s}(\bs u, \bs y)$ implies that,
for each $(\bs u, \bs y)\in D_r^1\times D_d^M\times D_s^N$, the 
function $P({\bs u}, {\bs y})$ has a convergent series expansion
about $(\bs u, \bs y)$.\\
\textit{iii)} Assumption \ref{ass} {\it ii)} implies that the dominant singularity of 
$G({\bs u}_0)$ is a square root. This type of singularity is generic for
functional equations with positive coefficients. If $P({\bs u},
{\bs y})$ is a \emph{polynomial} in $\bs u$ and $\bs y$, it
can be shown that Assumption \ref{ass} {\it ii)} follows from Assumption \ref{ass} 
{\it i)}, if $Q(u,y)$ is of degree $\ge2$ in $y$ and if $Q(u,0)\not\equiv0$.
By the closure properties of algebraic functions, it then follows, with an adaption
of Proposition \ref{prop:solv}, that all factorial moment generating functions are algebraic.\\
\textit{iv)} Assumption \ref{ass} {\it iii)} ensures that $G({\bs u_0})$ 
has exactly one singularity on its circle of convergence. The case of periodic
$G(\bs u_0)$ may be treated by a slight extension of this setup.

\medskip

We investigate analytic properties of $g_{\bf 0}(u_0)$. A function 
$f(u)$ is called \emph{$\Delta$-regular} \cite{FFK04} if it is analytic in the 
\emph{indented disc} $\Delta=\Delta(u_c,\eta,\phi)=\{u:|u|\le u_c+\eta, |\mbox{arg}
(u-u_c)|\ge\phi\}$ for some real numbers $u_c>0$, $\eta>0$ and $\phi$, 
where $0<\phi<\pi/2$. Note that $u_c\notin\Delta$, where we employ the 
convention $\mbox{arg}(0)=0$. The set of $\Delta$-regular functions is 
closed under addition, multiplication, differentiation, and integration. Moreover, 
if $f(u)\ne0$ in $\Delta$, then $1/f(u)$ exists in $\Delta$ and is $\Delta$-regular. The following
proposition is a straightforward extension of a well-known result, see 
e.g.~\cite[Thm.~10.6]{O95}, \cite[Prop.~1]{D97}, or \cite[Sec.~7.4]{FS07}.
 
\begin{prop}[cf.~\cite{O95,D97,FS07}]\label{prop:g0}
Given Assumption \ref{ass}, the power series $g_{\bf 0}(u_0)$ is analytic at $u_0=0$, 
with radius of convergence $0<u_c<\infty$. Its analytic continuation is $\Delta$-regular, 
with a square root singularity at $u_0=u_c$, and a local Puiseux expansion
\begin{equation}\label{form:genex0}
g_{\bf 0}(u_0) = g_{\bf 0}(u_c) + \sum_{l=0}^\infty f_{{\bf 0},l}
(u_c-u_0)^{-\gamma_{\bs0}+l/2},
\end{equation}
where $\gamma_{\bs0}=-1/2$, $f_{\bs 0,0}=-\sqrt{C/B}$, and $g_{\bf 0}(u_c)=
\lim_{u_0\to u_c^-} g_{\bs0}(u_0)<\infty$.
\qed
\end{prop}

\noindent \textbf{Remark.} The coefficients $f_{\bs 0,l}$, also called \emph{amplitudes},
can be computed recursively from the functional equation eqn.~\eqref{pgfeqn}, by inserting 
the representation eqn.~\eqref{form:genex0} into
eqn.~\eqref{pgfeqn} and then expanding the resulting equation in $s=\sqrt{u_c-u_0}$.
This technique will be exploited below, when using the method of dominant balance.

\medskip

The properties of $G({\bs u_0})$ carry over to the factorial moment generating functions 
$g_{\bs\nu}(u_0)$. We will first analyse the general form of the factorial moment generating 
functions, and later provide explicit values for exponents and leading amplitudes.

\begin{prop}\label{prop:pui}
Let Assumption \ref{ass} be satisfied. For $\bs\nu\neq\bs0$ arbitrary, the factorial moment 
generating function $g_{\bs\nu}(u_0)$ is analytic at $u_0=0$, 
with radius of convergence $0<u_c<\infty$. Its analytic continuation is $\Delta$-regular.
It has a local Puiseux expansion
\begin{equation}\label{form:genex}
g_{\bs\nu}(u_0) = \sum_{l=0}^\infty f_{{\bs\nu},l}(u_c-u_0)^{-\gamma_{\bs\nu}+l/2},
\end{equation}
with non-vanishing leading amplitude $f_{{\bs\nu},0}\ne0$ and exponent 
$\gamma_{\bs\nu}\in \frac{1}{2}\mathbb Z$.
\end{prop}

\begin{proof}
We prove the proposition by induction on $\bs\nu$ w.r.t.~the total order 
$\prec$ in Definition~\ref{def:totord}. Note that the proof of Proposition 
\ref{prop:solv} yields
\begin{equation}\label{form:gl1}
g_{\bs\nu}(u_0) Q_1(u_0,g_{\bf 0}(u_0))=Q_{\bs\nu}( 
u_0, g_{\bf 0}(u_0),\{ g_{\bs\mu}(u_0)\}_{\bs 0 \ne\bs\mu \prec\bs\nu}),
\end{equation}
where $Q_1(u,y)=1-\frac{\partial Q}{\partial y}(u,y)$ is analytic for 
$|u|<r$, $|y|<s$. The function $Q_{\bs\nu}( u,y,\widetilde 
{\bs y})$ is a polynomial in $\widetilde {\bs y}$ and 
analytic for $|u|<r$, $|y|<s$. Due to the closure properties of $\Delta$-regular and 
analytic functions, both $Q_1(u_0,g_{\bf 0}(u_0))$ and $Q_{\bs\nu}( 
u_0,g_{\bf 0}(u_0), \{ g_{\bs\mu}(u_0)\}_{\bs 0 \ne\bs\mu 
\prec\bs\nu})$ are $\Delta$-regular. Due to Proposition \ref{prop:g0}, we 
have $Q_1(u_0, g_{\bf 0}(u_0))\ne0$ in $\Delta$, thus its inverse is $\Delta$-regular. 
This implies that $g_{\bs\nu}(u_0)$ is $\Delta$-regular. 

Due to Proposition \ref{prop:g0}, the series $h(s)=y_c+\sum_{l=0}^\infty 
f_{{\bf 0},l}s^{(l+1)}$ is holomorphic at $s=0$ and equals $g_{\bf 0}(u_0)$ 
about $u_0=u_c$, if $s=\sqrt{u_c-u_0}$. We show that eqn.~\eqref{form:gl1} 
leads to local expansions of $g_{\bs\nu}(u_0)$ in terms of 
$s=\sqrt{u_c-u_0}$, which are meromorphic at $s=0$. Consider first the term 
$Q_1(u_0,g_{\bf 0}(u_0))$. The function $Q_1(u,y)$ is holomorphic at $(u,y)=
(u_c,y_c)$. Thus, as composition of holomorphic functions, $\widetilde 
Q_1(s):=Q_1(u_c-s^2,h(s))$ is holomorphic at $s=0$. We have $\widetilde 
Q_1(0)=0$ and $\widetilde Q'_1(0)=-\frac{\partial Q_1}{\partial y}(u_c,y_c) h'(0)
\ne0$. Thus, its inverse  $1/\widetilde Q_1(s)$ is meromorphic at $s=0$, 
with a simple pole. Consider finally the function $Q_{\bs\nu}
(u,y,\widetilde {\bs y})$. It is a polynomial in $\widetilde 
{\bs y}$ and analytic at $(u,y)=(u_c,y_c)$. Thus, after inserting 
the expansions of $g_{\bs\mu}(u_0)$, where we use Proposition 
\ref{prop:g0} and the induction hypothesis, we conclude that 
it is meromorphic in $s$ at $s=0$. 
It follows that $g_{\bs\nu}(u_0)$ is meromorphic in $s$ at $s=0$. 
Since $G(\bs u)$ is non-negative and $g_{\bs0}(u_0)$ 
does not vanish identically, $g_{\bs\nu}(u_0)$ does not vanish identically.  
We thus have local expansions \eqref{form:genex} of $g_{\bs\nu}(u_0)$.
\end{proof}

\noindent {\bf Remarks.} \textit{i)} 
It was argued in the preceding proof that $1-\sum_{j=1}^{N} \frac{\partial P}{\partial y_j}
(\bs u_0,{\bs G(\bs u_0)})=\widetilde Q_1(s)$ satisfies $\widetilde Q_1(s)\sim A
\sqrt{u_c-u_0}$ as $u_0\to u_c^{-}$, for some constant $A\ne0$, see also 
eqn.~\eqref{form:pref}. This will be used in the proof of the following proposition. \\
\textit{ii)} The exponents $\gamma_{\bs\nu}$ and, in principle,
all amplitudes $f_{\bs \nu,l}$ in the Puiseux expansion \eqref{form:genex} of 
$g_{\bs\nu}(u_0)$ can be computed recursively from the functional equation 
eqn.~\eqref{pgfeqn}, compare the argument above for the case $\bs\nu=0$. 
Our aim in this section is to determine the exponent $\gamma_{\bs\nu}$ and 
the amplitudes $f_{{\bs\nu},0}$. The method of Section \ref{sec:dom} will allow to 
also obtain information about the amplitudes $f_{{\bs\nu},l}$ for higher values of $l$,
with reasonable effort.

\medskip

We will first give an estimate of $\gamma_{\bs\nu}$.

\begin{prop}\label{prop:est}
Let Assumption \ref{ass} be satisfied. The exponent $\gamma_{\bs\nu}$ in 
the Puiseux expansion of $g_{\bs\nu}(u_0)$ eqns.~\eqref{form:genex0}, \eqref{form:genex},
satisfies the estimate
\begin{equation}
\gamma_{\bs\nu}\le -\frac{1}{2}+\sum_{i=0}^M \left( 1+\frac{i}{2}\right)\nu_i.
\end{equation}
\end{prop}

\noindent {\bf Remark.} Under mild additional assumptions, the above estimate
is sharp, i.e., we have $\gamma_{\bs\nu}= -\frac{1}{2}+\sum_{i=0}^M 
\left( 1+\frac{i}{2}\right)\nu_i$. This follows from Proposition \ref{prop:rek}, 
where we show that the amplitudes $f_{{\bs\nu},0}$, for the special choice of 
$\gamma_{\bs\nu}$ as above, have non-zero values. 
Similar estimates can be obtained in situations different from a square root 
singularity. The proof below can be adapted to treat these situations.

\begin{proof}
Set ${\widetilde\gamma_{\bs\nu}}
=-r+\sum_{i=0}^M r_i\nu_i$, where $r=1/2$, and $r_i=(1+i/2)$ for $i=0,\ldots,M$.
We have $r_{k+1}>r_k$ for $k=0,\ldots,M-1$. For 
$\bs\nu=\bf0$ and $\bs\nu={\bs e}_0$, we get 
from Proposition \ref{prop:g0} that $\gamma_{\bs\nu}=
{\widetilde\gamma_{\bs\nu}}$.
We prove the proposition by induction on $\bs\nu$ w.r.t.~the total order 
$\prec$ of Definition~\ref{def:totord}, using Proposition \ref{prop:pui}.   
For the remaining induction step, it suffices to show that $\gamma_{\bs\nu}
\le{\widetilde\gamma_{\bs\nu}}$.

\smallskip

Assume that $\gamma_{\bs\mu }\le{\widetilde\gamma_{\bs\mu }}$ 
holds for $\bs\mu \prec\bs\nu$. We group all terms in 
eqn.~\eqref{form:faa}, evaluated at $\bs u=\bs u_0$, 
which contain $g_{\bs\nu}(u_0)$, to the l.h.s. In the
resulting equation, its l.h.s.~$L(s)$ satisfies asymptotically, 
with $c_L$ a non-zero constant, $L(s)\sim c_L(u_c-u_0)^{-\gamma_{\bs
\nu}+1/2}$ as $u_0\to u_c^-$, see the remark following the proof of 
Proposition \ref{prop:pui}. By the reasoning in the proof of Proposition 
\ref{prop:pui}, the r.h.s.~$R(s)$ of the equation satisfies asymptotically, 
with $c_R$ a nonzero constant, $R(s)\sim c_R(u_c-u_0)^{-\gamma}$ as $u_0\to 
u_c^-$. Using the induction assumption, the exponent $\gamma$ can be estimated by
\begin{equation*}
\gamma \le \sum_{i=1}^s \gamma_{\bs\mu _i} |\bs\kappa_i|
\le \sum_{i=1}^s \left( -\frac{1}{2} + \sum_{j=0}^M r_j (\bs\mu _{i})_j\right)|
\bs\kappa_i|=  -\frac{|\bs\lambda|}{2}+\sum_{j=0}^M r_j\nu_j = 
- \frac{(|\bs\lambda|-1)}{2}+\widetilde\gamma_{\bs\nu}.
\end{equation*}
In the above expression, we used Lemma \ref{lemma:H} and the properties of 
$p_s(\bs\nu,\bs\lambda)$ in Lemma \ref{lemma:faa}.  Since 
$L(s)$ equals $R(s)$, we have $\gamma_{
\bs\nu}-1/2=\gamma$. Thus, the estimate $\gamma_{\bs\nu}\le
{\widetilde\gamma_{\bs\nu}}$ is satisfied for $|\bs\lambda |\ge2$. 
Let us analyse the terms contributing to $|\bs\lambda|=1$ in Lemma 
\ref{lemma:faa}. We have $s=1$, $|\bs\kappa_1|=1$ and $\bs\mu _1
=\bs\nu$. Now, use Lemma \ref{lemma:H}. If in eqn.~\eqref{form:chain} 
we have $|\bs\mu |=|\bs\nu|-l$ for some $l\in\{1,\ldots, |\bs\nu|\}$, 
it follows that $\widetilde\gamma_{\bs\mu }\le \widetilde\gamma_{\bs\nu}
-l\min_{0\le k\le M}\{r_k\}=\widetilde\gamma_{\bs\nu}-l$. Thus, terms with 
exponent larger than $\widetilde\gamma_{\bs\nu}-1/2$ must have $l=0$. 
If $|\bs\mu |=|\bs\nu|$ but $\bs\mu \ne \bs\nu$, 
we have $\widetilde\gamma_{\bs\mu }\le \widetilde\gamma_{\bs\nu}-
|r_j-r_k|$ for some indices $j, k\in\{0,\ldots, M\}$ satisfying $j\ne k$. 
This means that the exponent estimate $\gamma_{\bs\nu}-1/2\le-1/2+
\widetilde\gamma_{\bs\nu}$ is 
satisfied as long as $|r_j-r_k|\ge 1/2$ for $j, k\in\{0,\ldots,M\}$ and $j\ne k$. This 
condition is satisfied for the particular choice of exponents $r_i=(1+i/2)$ of 
the proposition. Thus, the proposition has been proved.
\end{proof}

The reasoning in the proof of Proposition \ref{prop:est} can be refined in 
order to obtain a recursion for the amplitudes $f_{{\bs\nu},0}$. 
For vectors $\bs x=(x_1,\ldots, x_K)\in\mathbb R^K$ and
$\bs y=(y_1,\ldots,y_K)\in\mathbb R^K$, we write 
${\bs x}\le {\bs y}$ if $x_i\le y_i$
for $i=1,\ldots,K$. In the proof, we will use the ``large oh'' symbol:
For functions $f(x)$ and $h(x)$, we write $f(x)={\cal O}(h(x))$ as 
$x\to x_0$ in a domain $D$, if there exists a positive constant $C$ and 
a neighbourhood $N(x_0)$ of $x_0$ such that $|f(x)|\le C |h(x)|$ for 
all $x\in N(x_0)\cap D$.  

\begin{prop}\label{prop:rek}
Let Assumption \ref{ass} be satisfied and fix $\bs\nu\in\mathbb N_0^{1+M}$.
Then, the amplitudes $f_{\bs\nu,0}=f_{\bs\nu}$ of the
Puiseux expansion eqn.~\eqref{form:genex} are, if ${\bs\nu}\ne {\bf 0}$ and 
${\bs\nu}\ne {\bs e}_0$, determined by the recursion
\begin{equation}\label{form:coefrek}
f_{\bs\nu} = \mu_0 \gamma_{{\bs\nu}-{\bs e}_1}
f_{{\bs\nu}-{\bs e}_1}+\sum_{i=1}^{M-1}\mu_i(\nu_i+1)
f_{{\bs\nu}-{\bs e}_{i+1}+{\bs e}_i}-\frac{1}
{2f_{\bf 0}}\sum_{\substack{{\bs\rho}\ne {\bf 0},{\bs\rho}
\ne {\bs\nu}\\ {\bf 0}\le{\bs\rho}\le{\bs\nu}}}
f_{\bs\rho}f_{{\bs\nu}-{\bs\rho}},
\end{equation}
with boundary conditions $f_{\bf 0}=-\sqrt{C/B}<0$, $f_{{\bs e}_0}
=-f_{\bf 0}/2>0$, and $f_{\bs\nu}=0$ if $\nu_j<0$ for some 
$j\in\{1,\ldots,M\}$. We have 
\begin{equation*}
\gamma_{\bs\nu}=-\frac{1}{2}+\sum_{i=0}^M \left(1+\frac{i}{2}\right)\nu_i.
\end{equation*}
For $i\in\{0,\ldots,M-1\}$, we have $\mu_i=-A_{i}/
(2Bf_{\bf 0})\ge0$, and the number $A_i\ge0$ is given by
\begin{equation}\label{eqn:Ai}
A_i=\sum_{j=1}^{N} \frac{\partial v_i^{(j)}}{\partial u_{i+1}}({\bs u}_c)
\frac{\partial P}{\partial y_j}({\bs u_c, \bs G(\bs u_c)}),
\end{equation}
where $\bs u_c=(u_c,1,\ldots,1)$. If $\bs\nu\ne {\bf 0}$, the amplitudes 
satisfy $f_{\bs\nu}\ge0$. This inequality is strict, if $A_i>0$ for $i=0,\ldots,M-1$.
\end{prop}

\begin{proof}
The amplitude $f_{\bs0}$ has been determined in Proposition \ref{prop:g0}.
Assume that ${\bs\nu}\ne {\bf 0}$ and ${\bs\nu}\ne {\bs e}_0$.
We group all terms in eqn.~\eqref{form:faa}, evaluated at $\bs u=\bs u_0$,
which contain $g_{\bs\nu}(u_0)$, to the l.h.s.
Using eqns.~\eqref{form:pref} and \eqref{form:facarmom}, the l.h.s.~$L(\bs u_0)$
of the resulting equation is given by
\begin{equation*}
L(\bs u_0)=\bs\nu! g_{\bs\nu}(u_0)P_{\bs e_1}({\bs u_0, \bs G(\bs u_0)})/
\left(\frac{{\rm d}G}{{\rm d} u_0}({\bs u}_0)\right).
\end{equation*}
Here, we used the identity
\begin{equation*}
\frac{{\rm d}G}{{\rm d}u_0}({\bs u}_0)=
\left(\sum_{j=1}^{N} \frac{\partial P}{\partial y_j}({\bs u_0, \bs G(\bs u_0)})
\right)\frac{{\rm d}G}{{\rm d}u_0}({\bs u}_0)+
P_{\bs e_1}({\bs u_0, \bs G(\bs u_0)}),
\end{equation*}
which is obtained from differentiating eqn.~\eqref{pgfeqn}.
By the reasoning in the proof of Proposition 
\ref{prop:est}, the r.h.s.~$R(u_0)$ of the resulting equation satisfies asymptotically, 
with $c_R$ a nonzero constant, $R(u_0)\sim c_R(u_c-u_0)^{-(\gamma_{\bs\nu}-1/2)}$ as $u_0\to 
u_c^-$. We collect all terms with exponents $\gamma_{\bs\nu}-1/2$. Due to 
the proof of Proposition \ref{prop:est}, they arise from terms where $|\bs\lambda|=1$ 
or $|\bs\lambda|=2$ in eqn.~\eqref{form:faa}. An explicit analysis using Lemma 
\ref{lemma:H}, whose details we omit, leads to the following expressions. The 
contribution arising from terms where $|\bs\lambda|=2$ is given by
\begin{equation*}
\frac{1}{2}\left(\sum_{j,k=1}^{N} \frac{\partial^2P}{\partial y_j\partial y_k}({\bs u_0, \bs 
G(\bs u_0)})\right)\bs\nu!
\sum_{\substack{\bs\rho\ne {\bf 0}, \bs\rho\ne \bs\nu\\{\bf 0}\le 
\bs\rho\le \bs\nu}}g_{\bs\rho}(u_0)g_{\bs\nu-\bs\rho}
(u_0).
\end{equation*}
The contribution from terms where $|\bs\lambda|=1$ is given by
\begin{equation*}
\sum_{j=1}^{N} \frac{\partial P}{\partial y_j}({\bs u_0, \bs G(\bs u_0)})
\sum_{i=0}^{M-1} (\bs\nu+{\bs e}_i)!
\frac{\partial v_i^{(j)}}{\partial u_{i+1}}({\bs u}_0)
g_{\bs\nu+{\bs e}_i-{\bs e}_{i+1}}(u_0).
\end{equation*}
Omitting arguments and normalising the l.h.s., we arrive at the equation
\begin{equation*}
\begin{split}
g_{\bs\nu}=&\frac{\sum \frac{\partial^2P}{\partial y_j\partial y_k}}{2P_{\bs e_1}}G'
\left( \sum g_{\bs\rho}g_{\bs\nu-\bs\rho}\right)
+\frac{G'}{P_{\bs e_1}}\sum_{j=1}^N \frac{\partial P}{\partial y_j}
\frac{\partial v_0^{(j)}}{\partial u_1}g_{\bs\nu-{\bs e}_1}'\\
&+\frac{G'}{P_{\bs e_1}}\sum_{j=1}^N \frac{\partial P}{\partial y_j} \sum_{i=1}^M (\nu_i+1)
\frac{\partial v_i^{(j)}}{\partial u_{i+1}} g_{\bs\nu+{\bs e}_i-{\bs e}_{i+1}}
+{\cal O}\left((u_c-u_0)^{-\gamma_{\bs\nu}+1/2}\right)
\end{split}
\end{equation*}
as $u_0\to u_c^-$, where $()'$ denotes differentiation w.r.t. $u_0$. We have $G'=g_{{\bs e}_1}$ 
and $f_{{\bs e}_1}=-f_{\bf 0}/2>0$. The implied recursion for the amplitudes $f_{{\bs\nu},0}$
is given by
\begin{equation*}
\begin{split}
f_{\bs\nu}=&\frac{\sum  \frac{\partial^2P}{\partial y_j\partial y_k}}{2P_{\bs e_1}}
\left(-\frac{f_{{\bf 0}}}{2}\right)\left( \sum f_{\bs\rho}f_{\bs\nu-\bs\rho}\right)
-\frac{f_{\bf 0}}{2P_{\bs e_1}}\sum_{j=1}^N \frac{\partial P}{\partial y_j}
\frac{\partial v_0^{(j)}}{\partial u_1}\gamma_{\bs\nu-{\bs e}_1} f_{\bs\nu-{\bs e}_1}\\
&-\frac{f_{{\bf 0}}}{2P_{\bs e_1}}\sum_{j=1}^N \frac{\partial P}{\partial y_j}
\sum_{i=1}^M (\nu_i+1)\frac{\partial v_i^{(j)}}{\partial u_{i+1}} f_{\bs\nu+{\bs e}_i-{\bs e}_{i+1}}.
\end{split}
\end{equation*}
After rewriting prefactors, using eqn.~\eqref{form:lett} and Proposition 
\ref{prop:pui}, we arrive at eqn.~\eqref{form:coefrek}.

We show strict positivity of the amplitudes, if $A_i>0$ for $i=0,\ldots,M-1$. 
For $|{\bs\nu}|\ge2$, all prefactors of $f_{\bs\mu }$ in eqn.~\eqref{form:coefrek} are
non-negative. For the first term, this follows from the estimate
\begin{equation*}
\gamma_{{\bs\nu}-{\bs e}_1}
= \left(-\frac{1}{2}-\frac{3}{2}+\sum_{i=0}^M\frac{i+2}{2} \nu_i \right)
\ge -2+|{\bs\nu}|\ge0.
\end{equation*}
Moreover, the last sum in eqn.~\eqref{form:coefrek} is over a nonempty set of 
indices $\bs\rho$ and has a strictly positive prefactor. For $|{\bs\nu}|=1$ and 
${\bs\nu}\ne {\bs e}_1$, the first term in eqn.~\eqref{form:coefrek} is zero, 
whereas the prefactors of $f_{\bs\mu }$ in the second term are strictly positive. 
Note finally that $f_{{\bs e}_1}=\mu_0\gamma_{\bf0}f_{\bf0}>0$, as is readily inferred from 
eqn.~\eqref{form:coefrek}. This leads, by induction, to strictly positive 
amplitudes $f_{\bs\nu}$. The same argument shows that $f_{\bs\nu}\ge{\bf 0}$
for $\bs\nu\ne\bf0$.
\end{proof}

\section{Method of dominant balance}\label{sec:dom}

We now discuss an altervative method to obtain the recursion for the 
amplitudes $f_{\bs\nu}$ in eqn.~\eqref{form:coefrek}, if 
$\nu_0=0$. It is based on a generating function approach, and generally 
easier to apply than the combinatorial approach in the proof of 
Proposition \ref{prop:rek}, which was based on an application of Faa 
di Bruno's formula. It also allows to analyse corrections to the 
leading singular behaviour of the factorial moment generating functions. 
Introduce parameters $\delta_i=1-u_i$, where $i\in\{1,\ldots,M\}$, and 
set $\bs\delta =(\delta_1, \ldots, \delta_M)$. From now on, we 
consider factorial moment generating functions $g_{\bs\nu}(u_0)$ 
for $\nu_0=0$ only. For ${\bs k}=(k_1,\ldots, k_M)\in\mathbb 
N_0^M$, set $g_{\bs k}(u_0)=g_{(0, \bs k)}(u_0)$, 
$f_{\bs k}(u_0)=f_{(0,\bs k)}(u_0)$, and $\gamma_{
\bs k}(u_0)=\gamma_{(0,\bs k)}(u_0)$.

\begin{prop}\label{prop:db1}
Assume that the $q$-functional equation~\eqref{form:funceq} has the solution 
$G(u_0,\bs u_+)\in{\cal H}_{d(P)}(\bs u_+)[[u_0]]$ such that $G(0,\bs u_+)=0$. 
Then ${\widetilde G}(u_0, \bs\delta):=G(u_0,1-\delta_1, \ldots, 1-\delta_M)\in
\mathbb C[[u_0,\bs\delta]]$ is a formal power series, given by
\begin{equation*}
{\widetilde G}(u_0,\bs\delta) = 
\sum_{{\bs k} \ge \bs 0} (-1)^{\bs k} 
g_{\bs k}(u_0) {\bs\delta}^{\bs k},
\end{equation*}
where $g_{\bs k}(u_0)=g_{(0, \bs k)}(u_0)\in\mathbb C[[u_0]]$ 
are the factorial moment generating functions \eqref{form:facarmom}.
\end{prop}

\begin{proof}
Proposition \ref{prop:uni} states that $G({\bs u})\in\mathcal{H}_{d(P)}(\bs u_+)[[u_0]]$.
Thus ${\widetilde G}(u_0,\bs\delta)\in\mathcal{H}_{d(P)-1}(\bs \delta)[[u_0]]\subset\mathbb 
C[[\bs\delta]][[u_0]]=\mathbb C[[u_0,\bs\delta]]$. We thus have
\begin{equation*}
{\widetilde G}(u_0,\bs\delta) = 
\sum_{\bs k\ge\bf 0} 
h_{\bs k}(u_0) {\bs\delta}^{\bs k},
\end{equation*}
for some $h_{\bs k}(u_0)\in\mathbb C[[u_0]]$. By Taylor's formula, we have
\begin{equation*}
h_{\bs k}(u_0)=\left.\frac{1}{\bs k!}
{\widetilde G}_{(0,\bs k)}
(u_0,\bs\delta)\right|_{\bs\delta={\bf  0}}=
\left.\frac{(-1)^{\bs k}}{\bs k!}
G_{(0,\bs k)}({\bs u})\right|_{\bs u= \bs u_0}=
(-1)^{\bs k}
g_{\bs k}(u_0),
\end{equation*}
and the proposition is proved.
\end{proof}

Let $\mathbb C((s))$ denote the field of formal Laurent series $f(s)=
\sum_{l\ge l_0}f_l s^l$, where $l_0\in\mathbb Z$ and $f_l\in\mathbb C$ 
for $l\ge l_0$. Employing the Puiseux expansions of the factorial moment 
generating functions, we have an alternative representation of the 
generating function $G({\bs u})$.

\begin{prop}\label{prop:scalbeh}
Let Assumption \ref{ass} be satisfied.
Replace the coefficients $g_{\bs k}(u_0)$ of ${\widetilde G}
(u_0, \bs\delta)\in\mathbb C[[u_0,\bs\delta]]$ in 
Proposition \ref{prop:db1} by their Puiseux expansions of Propositions 
\ref{prop:g0} and \ref{prop:pui}, and denote the resulting series by 
${\widetilde G}(s, \bs\delta)\in\mathbb C((s))[[\bs \delta]]$.
It is explicitly given by
\begin{equation*}
{\widetilde G}(s, \bs\delta)=G({\bs u_c})+
\sum_{\bs k\ge\bf0}\left(\sum_{l=0}^\infty
(-1)^{\bs k}f_{{\bs k},l}s^{-\gamma_{\bs k}+l/2}
\right) \bs\delta^{\bs k},
\end{equation*}
where $\gamma_{\bs k}=-1/2+\sum_{i=1}^M(1+i/2)k_i$,
where the numbers $f_{{\bs k},l}=f_{(0,{\bs k}),l}$ are 
defined in eqn.~\eqref{form:genex0} and in eqn.~\eqref{form:genex},
and where $\bs u_c=(u_c,1,\ldots,1)$.

\smallskip

Then, the rescaled series $G({\bs u}(s,\bs\epsilon))=
{\widetilde G}(s,\epsilon_1s^3,\epsilon_2s^4,\ldots,\epsilon_Ms^{M+2})
\in\mathbb C[[s,\bs \epsilon]]$ is a \emph{formal power series},
\begin{equation}\label{form:scaling}
G({\bs u}(s,\bs\epsilon)) = G({\bs u_c}) + s F\left(s,
\bs\epsilon\right),
\end{equation}
where $F(s,\bs\epsilon)\in\mathbb C[[s, \bs\epsilon]]$ is 
given by
\begin{equation}\label{form:Fl}
F(s, \bs\epsilon)=\sum_{l=0}^\infty F_l({\bs\epsilon}) s^l,
\qquad F_l(\bs\epsilon)=\sum_{\bs k\ge\bf0} 
(-1)^{\bs k}f_{{\bs k},l}
{\bs\epsilon}^{\bs k}.
\end{equation}
\end{prop}

\begin{proof}
Due to Propositions \ref{prop:g0} and \ref{prop:pui}, we have 
$g_{\bs k}(u_0)\in\mathbb C((s))$, where $s=\sqrt{u_c-u_0}$. The 
explicit form \eqref{form:scaling} follows immediately from the Puiseux expansions 
in Propositions \ref{prop:g0} and \ref{prop:pui}, together with
the exponent estimate in Proposition \ref{prop:est}. 
\end{proof}

\noindent {\bf Remarks.}
$i)$ Equations like eqn.~\eqref{form:scaling} appear in statistical mechanics as
a so-called \emph{scaling Ansatz}, being an assumption on the behaviour of a generating
function near a multicritical point singularity \cite{J00}, see also \cite{BFG03,KS05}. Its validity 
has been proved only in a limited number of examples. Here, we employ the 
different framework of formal power series.\\
$ii)$ Proposition \ref{prop:scalbeh} suggests an alternative strategy to compute the 
singular behaviour of the factorial moment generating functions. We consider the
functional equation for $F(s,\bs\epsilon)\in\mathbb C[[s, \bs\epsilon]]$,
which is induced by the $q$-functional equation for $G(\bs u)$. Writing
$F(s, \bs\epsilon)=\sum_{l\ge0} F_l({\bs\epsilon}) s^l$,
where $F_l({\bs\epsilon}) \in\mathbb C[[\bs\epsilon]]$, then
leads to a partial differential equation for $F_l({\bs\epsilon})$, upon 
expanding the induced functional equation in powers of $s$. Note that this algorithm,
which is computationally involved, can be easily implemented in a computer algebra 
system. The above method is called in statistical mechanics the method of {\em dominant 
balance}, see \cite{PB95a,R02}. It can be applied in our framework, if an exponent 
bound like that of Proposition \ref{prop:est} is known. Such bounds are generally 
easier to obtain than explicit recursions for amplitudes, see the above proofs.

\begin{prop}\label{prop:newfeqn}
Let Assumption \ref{ass} be satisfied.
Define the power series ${\bs u}(s,\bs\epsilon)=(u_0(s,\bs\epsilon), 
\ldots, u_M(s,\bs\epsilon))$, where
\begin{equation*}
u_0(s,\bs\epsilon)=u_c-s^2, \qquad 
u_i(s,\bs\epsilon)=1-\epsilon_is^{i+2} \qquad (i=1,\ldots,M).
\end{equation*}
For a $q$-shift $\bs v$, define the induced $q$-shift $s_{{\bs v}}$ of $s$ and the 
induced $q$-shift $\bs\epsilon_{{\bs v}}=(\epsilon_{1,{\bs v}}, \ldots, 
\epsilon_{M,{\bs v}})$ of $\bs\epsilon$ by
\begin{equation*}
\begin{split}
s_{{\bs v}}&=s_{{\bs v}}(s,\bs\epsilon)=
\sqrt{u_c-v_0({\bs u}(s,\bs\epsilon))},\\
\epsilon_{i,{\bs v}}
&=\epsilon_{i,{\bs v}}(s,\bs\epsilon)=
\frac{1-v_i({\bs u}(s,\bs\epsilon))}{s_{{\bs v}}
(s,\bs\epsilon)^{i+2}} \qquad (i=1,\ldots,M).
\end{split}
\end{equation*}
We then have $s_{{\bs v}}\in\mathbb C[[s, \bs\epsilon]]$ and 
$\epsilon_{i,{\bs v}}\in\mathbb C[[s, \bs\epsilon]]$ for $i=1,\ldots,M$.
The functional equation for $F(s, \bs\epsilon)\in\mathbb C[[s, 
\bs\epsilon]]$ of Proposition \ref{prop:scalbeh}, induced by 
eqn.~\eqref{form:funceq}, is given by
\begin{equation}\label{form:funcind}
G({\bs u}(s,\bs\epsilon))=P( {\bs u}(s,\bs\epsilon),
{\bs H}({\bs u}(s,\bs\epsilon))).
\end{equation}
In the above equation, $G({\bs u}(s,\bs\epsilon))$ is given by 
eqn.~\eqref{form:scaling}, and the power series $H^{(j)}({\bs u}
(s,\bs\epsilon))$ are given by
\begin{equation*}
H^{(j)}({\bs u}(s,\bs\epsilon)) = G({\bs u_c})+s_{{\bs v}^{(j)}}
F(s_{{\bs v}^{(j)}},\bs\epsilon_{{\bs v}^{(j)}}) \qquad (1\le j\le N).
\end{equation*}
\end{prop}

\begin{proof}
This follows from direct computation. Note first that the Taylor expansions 
of the power series $v_k({\bs u})$ about ${\bs u}={\bs u}_0$ 
are of the form
\begin{equation}\label{form:taylvk}
v_k({\bs u})= v_k({\bs u}_0) +
\sum_{\substack{\bs l\ne 0\\ \bs l\ge 0}}
\frac{(-1)^{\bs l}}{\bs l !} v_{k,\bs l}
({\bs u}_0)(1-u_1)^{l_1}\cdot\ldots\cdot(1-u_M)^{l_M} \qquad (k=0,\ldots,M).
\end{equation}
Now insert the parametrisations $u_k(s,\bs\epsilon)\in \mathbb 
C[s,\bs\epsilon]$, where $k=0,\ldots,M$. We get $v_k({\bs u}
(s,\bs\epsilon))\in\mathbb C[[{s,\bs\epsilon}]]$. Due to the 
$q$-shift properties, we have
\begin{equation*}
\begin{split}
&v_0({\bs u}(s,\bs\epsilon))=u_c-s^2-\frac{\partial v_0}{\partial u_1}
({\bs u}_c) \epsilon_1s^3+{\cal O}(s^4),\\
&v_k({\bs u}(s,\bs\epsilon))= 1-\epsilon_ks^{k+2}-
\frac{\partial v_k}{\partial u_{k+1}}({\bs u}_c)\epsilon_{k+1}s^{k+3} + 
{\cal O}(s^{k+4}) \qquad (k=1,\ldots,M-1),\\
&v_M({\bs u}(s,\bs\epsilon))= 1-\epsilon_Ms^{M+2}+ {\cal O}(s^{M+3}).
\end{split}
\end{equation*}
We thus have $v_0({\bs u}(s,\bs\epsilon))=u_c-s^2-s^3R_0(s,\bs
\epsilon)$ and $v_k({\bs u}(s,\bs\epsilon))=1-s^{k+2}R_k(s,\bs
\epsilon)$ for $k=1,\ldots,M$, where $R_i(s,\bs\epsilon)\in\mathbb 
C[[s, \bs\epsilon]]$ for $i=0,\ldots,M$. For $s_{{\bs v}}$, we get
\begin{equation*}
s_{{\bs v}}=\sqrt{u_c-v_0({\bs u}(s,\bs\epsilon))}=
s\sqrt{1-sR_0(s,\bs\epsilon)}.
\end{equation*}
We thus have $s_{{\bs v}}\in\mathbb C[[s, \bs\epsilon]]$ 
and $s_{{\bs v}}(0,\bs\epsilon)=0$. We get
\begin{equation*}
\epsilon_{i,{\bs v}}=\frac{1-v_i({\bs u}(s,\bs\epsilon))}
{s_{{\bs v}}^{i+2}}=\frac{R_i(s,\bs\epsilon)}{\sqrt{1-sR(s,
\bs\epsilon)}^{i+2}} \qquad (i=1,\ldots,M).
\end{equation*}
It follows that $\epsilon_{i,{\bs v}}\in\mathbb C[[s, \bs\epsilon]]$, 
where $i=1,\ldots,M$. Furthermore, we get from eqn.~\eqref{form:taylvk}
\begin{equation*}
\epsilon_{i,{\bs v}}(0,\bs\epsilon)=R_i(0,\bs\epsilon)=
\epsilon_i,
\end{equation*}
such that $\bs\epsilon_{{\bs v}}(0,\bs\epsilon)=\bs\epsilon$. 
This ensures that eqn.~\eqref{form:funcind} is well-defined for  
$F(s, \bs\epsilon)\in\mathbb C[[s, \bs\epsilon]]$.
\end{proof}

\begin{prop}\label{theo:pde}
Let Assumption \ref{ass} be satisfied. The formal power series
$F_0(\bs\epsilon)\in\mathbb C[[\bs\epsilon]]$
satisfies the singular, quasi-linear partial differential equation of first order
\begin{equation}\label{form:res}
\sum_{j=1}^{N} \frac{\partial P}{\partial y_j}(\bs u_c, \bs G(\bs u_c)) \left( \frac{1}{2}
\frac{\partial v_0^{(j)}}{\partial u_1}({\bs u}_c)\epsilon_1 F_0(\bs\epsilon) +
\sum_{i=1}^M h_i^{(j)}(\bs\epsilon) \frac{\partial}{\partial\epsilon_i}
F_0(\bs\epsilon) \right) +B F_0^2(\bs\epsilon)-C=0,
\end{equation}
where the numbers $B$, $C$ are defined in eqn.~\eqref{form:lett}, and the 
formal power series $h_i^{(j)}(\bs\epsilon)$ are given by
\begin{equation}\label{form:hi}
\begin{split}
&h_i^{(j)}(\bs\epsilon)=\frac{\partial v_i^{(j)}}{\partial u_{i+1}}({\bs u}_c)
\epsilon_{i+1}-\frac{i+2}{2}\frac{\partial v_0^{(j)}}{\partial u_{1}}({\bs u}_c)
\epsilon_1\epsilon_i \qquad (i=1,\ldots, M-1),\\
&h_M^{(j)}(\bs\epsilon)=-\frac{M+2}{2}\frac{\partial v_0^{(j)}}
{\partial u_{1}}({\bs u}_c)\epsilon_1\epsilon_M.
\end{split}
\end{equation}
\end{prop}

\begin{proof}
Using the expansions in the proof of Proposition \ref{prop:newfeqn}, it 
is readily inferred that the power series $s_{{\bs v}}$ and $\bs
\epsilon_{{\bs v}}$ are, to leading orders in $s$, given by
\begin{equation*}
\begin{split}
&s_{{\bs v}^{(j)}}= s+ \frac{1}{2}\frac{\partial v_0^{(j)}}
{\partial u_{1}}({\bs u}_c)\epsilon_1 s^2 +{\cal O}(s^3),\\
&\epsilon_{i,{\bs v}^{(j)}}= \epsilon_i+h_i^{(j)}(\bs\epsilon)s
+{\cal O}(s^2) \qquad (i=1,\ldots, M-1),
\end{split}
\end{equation*}
with $h_i^{(j)}(\bs\epsilon)$ as defined in eqn.~\eqref{form:hi}.
Using this result, we compute the expansion of $H^{(j)}({\bs u}(s,\bs\epsilon))$
up to order $s^2$. This yields
\begin{equation*}
\begin{split}
H^{(j)}(&{\bs u}(s,\bs\epsilon))  = G({\bs u_c})+
s_{{\bs v}^{(j)}}F(s_{{\bs v}^{(j)}},\bs
\epsilon_{{\bs v}^{(j)}})\\
&= G({\bs u_c}) + F_0(\bs\epsilon)s +
\left( F_1(\bs\epsilon)+\frac{1}{2}\frac{\partial v_0^{(j)}}
{\partial u_{1}}({\bs u}_c)\epsilon_1 F_0(\bs\epsilon)
+ \sum_{i=1}^M h_i^{(j)}(\bs\epsilon) \frac{\partial}{\partial\epsilon_i}
F_0(\bs\epsilon)\right)s^2 +{\cal O}(s^3).
\end{split}
\end{equation*}
Now expand the functional equation \eqref{form:funcind} to leading orders in 
$s$. Terms of order $s^0$ vanish due to eqn.~\eqref{pgfeqn}, evaluated at 
${\bs u_0}={\bs u_c}$. Terms of order $s^1$ vanish due the 
condition $\sum_{j=1}^{N} \frac{\partial P}{\partial y_j}(\bs u_c, \bs 
G(\bs u_c))=1$. Terms of order $s^2$ lead to the partial differential 
equation given above.
\end{proof}

\noindent {\bf Remarks.}
{\it i)} For $M>2$, it is an open question whether closed form solutions 
for $F_0(\bs\epsilon)$ exist. See \cite{NT03b} for a discussion of the cases
$M=1$ and $M=2$.\\
{\it ii)} The above method can also be used to derive 
partial differential equations characterising the generating functions 
$F_l(s,\bs\epsilon)$ of the amplitudes $f_{{\bs k},l}$ 
for $l>0$. These equations arise in the expansion of the $q$-functional 
equation in $s$ at order $l+2$, see \cite{R02} for examples where $M=1$.

\medskip

The above theorem leads to an alternative derivation of the recursion 
eqn.~\eqref{form:coefrek} of Proposition \ref{prop:rek} in the case 
$\nu_0=0$.

\begin{proof}[Alternative proof of eqn.~\eqref{form:coefrek}]
We set $F_0(\bs\epsilon)=K(-\bs\epsilon)+f_{{\bf 0},0}$ 
and rewrite eqn.~\eqref{form:res} in the form
\begin{equation}\label{form:ee1}
K(\bs\epsilon)=\sum_{i=1}^M \frac{i+2}{2}\mu_0\epsilon_1\epsilon_i
\frac{\partial K(\bs\epsilon)}{\partial \epsilon_i}+\sum_{i=1}^{M-1} 
\mu_i\epsilon_{i+1}\frac{\partial K(\bs\epsilon)}{\partial \epsilon_i}
-\frac{\mu_0}{2}\epsilon_1(K(\bs\epsilon)+f_{{\bf 0},0})-\frac{1}
{2f_{{\bf 0},0}}K(\bs\epsilon)^2,
\end{equation}
where $K({\bf 0})=0$, and the constants $\mu_i$ are, for $i\in\{0,\ldots,M-1\}$, 
given by $\mu_i=-A_i/(2Bf_{{\bf 0},0})$, with $A_i$ defined in eqn.~\eqref{eqn:Ai}.
This leads to the recursion eqn.~\eqref{form:coefrek} for the coeffients 
$f_{\bs\nu}$ in Proposition \ref{prop:rek}, if we set $\bs
\nu=(0,{\bs k})$.
\end{proof}

\section{Growth of amplitudes}\label{sec:growth}

We are interested in the growth of the coefficents $f_{\bs\nu,0}$,
which appear in Proposition \ref{prop:rek} and in Proposition \ref{theo:pde}, 
in the case $\nu_0=0$. To this end, we study properties of the partial 
differential equation of the associated generating function $F_0(\bs
\epsilon)$ of eqn.~\eqref{form:Fl}, given by
\begin{equation*}
F_0(\bs\epsilon)=\sum_{\bs k\ge\bf0} (-1)^{\bs k}
f_{{\bs k},0}{\bs\epsilon}^{\bs k}.
\end{equation*}

\begin{prop}\label{theo:gr}
For a $q$-functional equation, let Assumption \ref{ass} be satisfied.
There exist positive real numbers $D, R_1,\ldots,R_M$, such that
\begin{equation*}
|f_{{\bs k},0}| \le D (k_1+\ldots+k_M)! \, (R_1)^{k_1}\cdot\ldots
\cdot(R_M)^{k_M}
\end{equation*}
for all $\bs k\ge\bf0$.
\end{prop}

For the proof of the proposition, we apply the technique of majorising series.
A formal power series $g=\sum_{\bs k} g_{\bs k} 
{\bs x}^{\bs k}$ majorises a formal power series $h=
\sum_{\bs k} h_{\bs k} {\bs x}^{\bs k}$ if 
$|g_{\bs k}|\le |{\bs k}|!h_{\bs k}$ for all 
${\bs k}\ge{\bf 0}$. We then write $g\ll h$. We have the following relations.

\begin{lemma}\label{lem:maj}
Let $g\ll h$. Then 
\begin{equation}
g^2\ll h^2, \qquad x_1g\ll x_1h, \qquad
x_1x_i\frac{\partial g}{\partial x_i}\ll x_1 h, \qquad
x_{j+1}\frac{\partial g}{\partial x_j}\ll x_{j+1}\frac{\partial h}{\partial x_j},
\end{equation}
for $i=1,\ldots,M$ and $j=1,\ldots, M-1$.
\end{lemma}

\begin{proof}[Sketch of proof]
These relations are checked by direct computation. If $i=2,\ldots,M$, we have
$x_1x_i\frac{\partial g}{\partial x_i}=\sum k_ig_{{\bs k}-
{\bs e}_1}{\bs x}^{\bs k}$. Thus $|k_ig_{{\bs k}-
{\bs e}_1}|\le (|{\bs k}|-1)! k_ih_{{\bs k}-
{\bs e}_1}\le |{\bs k}|!h_{{\bs k}-{\bs e}_1}$, 
and the statement follows. The remaining assertations are proved similarly.
\end{proof}

\begin{proof}[Proof of Proposition \ref{theo:gr}]
Starting with $K(\bs\epsilon)$, defined in eqn.~\eqref{form:ee1}, 
we introduce the majorant equation
\begin{equation}\label{form:ee2}
L(\bs\epsilon)=\mu_0\left(\sum_{i=1}^M \frac{i+2}{2}\right) \epsilon_1 
L(\bs\epsilon)
+\frac{\mu_0}{2}\epsilon_1(L(\bs\epsilon)+|f_{{\bf 0},0}|)
+\frac{1}{2|f_{{\bf 0},0}|}
L(\bs\epsilon)^2
+ \sum_{i=1}^{M-1} \mu_i\epsilon_{i+1}
\frac{\partial L(\bs\epsilon)}{\partial \epsilon_i}.
\end{equation}
The last equation uniquely defines a power series with non-negative 
coefficients and positive radius of convergence, satisfying $L({\bf 0})=0$. 
It belongs to a class of singular partial differential equations with regular 
solutions, which is discussed in 
\cite[Thm.~2.8.2.1 and Sec.~2.9.5]{GT96}. To prove that 
$K(\bs\epsilon)\ll L(\bs\epsilon)$, denote the 
r.h.s. of eqn.~\eqref{form:ee1} by $RK$ and the r.h.s.~of eqn.~\eqref{form:ee2} by 
$\widehat RL$. By construction, $g\ll h$ implies $Rg \ll \widehat R h$. We use 
an iteration argument. Set $K_0=L_0=0$. Clearly $K_0\ll L_0$. Define 
$K_n=RK_{n-1}$ and $L_n=\widehat R L_{n-1}$ for $n\in\mathbb N$. We have 
$K_n\ll L_n$ for $n\in\mathbb N_0$, due to Lemma \ref{lem:maj}. Thus 
$K(\bs\epsilon)\ll L(\bs\epsilon)$ for the 
formal solutions $K(\bs\epsilon)$ of eqn.~\eqref{form:ee1} and 
$L(\bs\epsilon)$ of eqn.~\eqref{form:ee2}. Since 
$L(\bs\epsilon)$ is regular, we have the estimate
\begin{equation*}
|f_{{\bs k},0}| \le |{\bs k}|! [\bs\epsilon^{\bs k}]
L(\bs\epsilon)\le D |{\bs k}|! {\bs R}^{\bs k}
\end{equation*}
for some real constants $D>0$ and $R_i>0$, where $i=1,\ldots,M$.
\end{proof}

\section{Moments and limit distributions}\label{sec:mom}

In the following, we give a probabilistic interpretation of the obtained results.
This generalises the discussion of Dyck paths in the introduction, before 
Proposition~\ref{prop:dyckmom}, to the case of a general $q$-functional equation,
and will prove part \textit{i)} of Theorem~\ref{theo:probdist}. For technical reasons 
(Lemma~\ref{lem:Levy}), we will not use random variables below, but rather argue
with the associated probability measures.

\smallskip

For a $q$-functional equation, let Assumption \ref{ass} be satisfied. Then, 
the coefficients $p_{\bs n}$ in a solution $G({\bs u})= \sum_{\bs 
n\ge 0} p_{\bs n}{\bs u^n}$, such that $G(\bs0)=0$, are non-negative. Assume that the 
numbers $A_i$ of Proposition \ref{prop:rek} satisfy $A_i>0$ for $i=0,\ldots, M-1$. 
We then have $f_{{\bs k},0}>0$ for all ${\bs k}\ge {\bf 0}$ and 
${\bs k}\neq {\bf 0}$. This implies that $\sum_{n_1,\ldots,n_M}
p_{n_0,n_1,\ldots,n_M}>0$ for almost all $n_0$, see eqn.~\eqref{form:coefas} below. 
Fix $N_0\in\mathbb N$ such that strict positivity holds for all $n_0\ge N_0$. 
For $n_0\ge N_0$, we define discrete Borel probability measures $\widetilde\mu_{n_0}$ 
by
\begin{equation*}
\widetilde\mu_{n_0}=\sum_{n_1,\ldots,n_M} \frac{p_{n_0,n_1,\ldots,n_M}}{\sum_{m_1,
\ldots,m_M} p_{n_0,m_1,\ldots,m_M}} \delta_{(n_1,\ldots,n_M)} ,
\end{equation*}
compare eqn.~\eqref{eq:pxt} in the introduction.
Their corresponding moments are, for ${\bs k}\in\mathbb N_0^M$, given by
\begin{equation}\label{form:tilmom}
\widetilde m_{\bs k}(n_0) = \frac{\sum_{n_1,\ldots,n_M} n_1^{k_1}\cdot\ldots
\cdot n_M^{k_M} p_{n_0,n_1,\ldots,n_M}}{\sum_{n_1,\ldots,n_M} p_{n_0,n_1,\ldots,n_M}}.
\end{equation}

We are interested in the asymptotic behaviour of the moments $\widetilde 
m_{\bs k}(n_0)$, as $n_0$ tends to infinity. This will be obtained by an 
analysis of the coefficients of the factorial moment generating functions, 
which relies on a transfer lemma \cite[Thm.~1]{FO90}.

\begin{lemma}[Transfer lemma \cite{FO90}]
Suppose that $F(z)=\sum_{n\ge0}f_nz^n$ has a singularity at $z=z_c$ and is
$\Delta$-regular, i.e., it is analytic in the domain
\begin{equation*}
\Delta=\{z:|z|\le z_c+\eta, |{\rm arg}(z-z_c)|\ge\phi\}
\end{equation*}
for some $\eta>0$ and $0<\phi<\pi/2$. Assume that, as $z\to z_c$ in $\Delta$,
\begin{equation*}
F(z) = {\cal O} \left( (z_c-z)^\alpha \right)
\end{equation*}
for some real $\alpha$. Then, the $n$-th Taylor coefficient $f_n$ of $F(z)$ satisfies
\begin{equation*}
f_n={\cal O}(z_c^{-n} n^{-1-\alpha}) \qquad (n\to\infty).
\end{equation*}
\qed
\end{lemma}

The following lemma characterises the asymptotic behaviour
of the moments $\widetilde m_{\bs k}(n_0)$, as $n_0$ tends to infinity.

\begin{lemma}\label{lem:pos}
Let Assumption \ref{ass} be satisfied. Assume that the 
numbers $A_i$ of Proposition \ref{prop:rek} satisfy $A_i>0$ for 
$i=0,\ldots, M-1$. Then the  moments $\widetilde m_{\bs k}(n_0)$ 
eqn.~\eqref{form:tilmom} are well-defined for almost all $n_0$. They  
are for ${\bs k}\in\mathbb N_0^M$ asymptotically given by
\begin{equation*}
\widetilde m_{\bs k}(n_0)= \frac{{\bs k}!}{f_{{\bf 0},0}
u_c^{\gamma_{\bs k}-\gamma_{\bf 0}}}
\frac{\Gamma(\gamma_{\bf 0})}{\Gamma(\gamma_{\bs k})}f_{{\bs k},0}
n_0^{\gamma_{\bs k}-\gamma_{\bf 0}}+{\cal O}(n_0^{\gamma_{\bs k}-
\gamma_{\bf 0}-1/2}) \qquad (n_0\to\infty),
\end{equation*}
where $\Gamma(z)$ denotes the Gamma function, and where the numbers 
$f_{{\bs k},0}$ and $\gamma_{\bs k}$ are defined in 
Proposition \ref{prop:rek}.
\end{lemma}

\begin{proof}
The functions $g_{\bs k}(u_0)$ are $\Delta$-regular due to 
Proposition \ref{prop:pui}. We infer from eqn.~\eqref{form:genex} that 
$g_{\bs k}(u_0)=f_{{\bs k},0}(u_c-u_0)^{-\gamma_{\bs k}}+
{\cal O}((u_c-u_0)^{-\gamma_{\bs k}+1/2})$ as $u_0\to u_c^-$, where 
$\gamma_{\bf 0}=-1/2$ and $\gamma_{\bs k}>0$ otherwise. The 
coefficient asymptotics of $(u_c-u_0)^{-\gamma_{\bs k}}$ and 
the transfer lemma yield
\begin{equation}\label{form:coefas}
[u_0^{n_0}] g_{\bs k}(u_0)=\frac{f_{{\bs k},0}}
{u_c^{\gamma_{\bs k}}\Gamma(\gamma_{\bs k})}u_c^{-n_0}
n_0^{\gamma_{\bs k}-1}\left(1+{\cal O}(n_0^{-1/2})\right) \qquad
(n_0\to\infty),
\end{equation}
where $[x^n]f(x)$ denotes the coefficient of $x^n$ in the power 
series $f(x)$. The error term implies that asymptotically 
factorial moments coincide with ordinary moments. The numbers 
$\widetilde m_{\bs k}(n_0)$ are thus well-defined for $n_0$ 
large enough, and asymptotically given by
\begin{equation*}
\begin{split}
\widetilde m_{\bs k}(n_0)&=
\frac{\sum_{{\bs n}_+\ge0} {\bs n}_+^{{\bs k}} 
p_{n_0,{\bs n}_+}}{\sum_{{\bs n}_+\ge0} 
p_{n_0,{\bs n}_+}}=\frac{[u_0^{n_0}]g_{\bs k}(u_0)}
{[u_0^{n_0}]g_{\bf 0}(u_0)}{\bs k}!\left( 1+{\cal O}(n_0^{-1/2}) 
\right)\\
&= \frac{{\bs k}!}{f_{{\bf 0},0}u_c^{\gamma_{\bs k}-
\gamma_{\bf 0}}}\frac{\Gamma(\gamma_{\bf 0})}{\Gamma(\gamma_{\bs k})}
f_{{\bs k},0}n_0^{\gamma_{\bs k}-\gamma_{\bf 0}}\left( 1+
{\cal O}(n_0^{-1/2}) \right) \qquad (n_0\to\infty).
\end{split}
\end{equation*}
This concludes the proof of the lemma.
\end{proof}
The moments $\widetilde m_{\bs k}(n_0)$ diverge as $n_0\to\infty$,
as may be inferred from Lemma \ref{lem:pos}. Introduce normalised Borel
probability measures $\mu_{n_0}$ by
\begin{equation}\label{form:normrv}
\mu_{n_0}=\sum_{n_1,\ldots,n_M} \frac{p_{n_0,n_1,\ldots,n_M}}{\sum_{m_1,\ldots,m_M} 
p_{n_0,m_1,\ldots,m_M}} \delta_{(\widetilde n_1,\ldots,\widetilde n_M)},
\end{equation}
where $\widetilde n_k=n_kn_0^{-(k+2)/2}$ for $k=1,\ldots,M$. For ${\bs k}
\in\mathbb N_0^M$, denote the corresponding moments by $m_{\bs k}(n_0)$.
We now show that the limits
\begin{equation}\label{form:normmom}
m_{\bs k}=\lim_{n_0\to\infty} m_{\bs k}(n_0)=\lim_{n_0\to\infty} 
\frac{\widetilde m_{\bs k}(n_0)}{n_0^{\gamma_{\bs k}-\gamma_{\bf 0}}}
= \frac{{\bs k}!}{f_{{\bf 0},0}u_c^{\gamma_{\bs k}-\gamma_{\bf 0}}}
\frac{\Gamma(\gamma_{\bf 0})}{\Gamma(\gamma_{\bs k})}f_{{\bs k},0}
\end{equation}
exist and define a unique Borel probability measure $\mu$, with finite moments 
$m_{\bs k}>0$ at all orders. This will be achieved using L\'evy's 
continuity theorem. We first prove the following lemma.

\begin{lemma}\label{lem:momcon}
Let Assumption \ref{ass} be satisfied. Assume that the 
numbers $A_i$ of Proposition \ref{prop:rek} satisfy $A_i>0$ for $i=0,\ldots, M-1$. 
Consider for ${\bs k}\in \mathbb N_0^M$
the numbers $m_{\bs k}\ge0$, defined in eqn.~\eqref{form:normmom} and
eqn.~\eqref{form:tilmom}.
For ${\bs t}\in \mathbb R^M$, we have
\begin{equation*}
\lim_{|{\bs k}|\to\infty} \frac{m_{\bs k}
{\bs t}^{\bs k}}{{\bs k}!}=0.
\end{equation*}
Equivalently, the exponential generating function of the
numbers $m_{\bs k}$ is entire.
\end{lemma}

\begin{proof}
The limit $m_{\bs k}$ eqn.~\eqref{form:normmom} exists for ${\bs k}\in 
\mathbb N_0^M$, due to Lemma~\ref{lem:pos}. Note that
\begin{equation*}
\gamma_{\bs k}=-\frac{1}{2}+\sum_{i=1}^M \left(1+\frac{i}{2}\right)k_i\ge
-\frac{1}{2}+\frac{3}{2}|{\bs k}|\ge |{\bs k}|+1+ \lfloor \frac{|{\bs k}|-3}{2}\rfloor.
\end{equation*}
Thus $|{\bs k}|\to\infty$ implies $\gamma_{\bs k}\to\infty$. Furthermore, since
$e(n/e)^n\le n!\le e n (n/e)^n$ for $n\in\mathbb N$, we have the estimate
\begin{equation*}
\frac{n!}{(n+n_0)!}\le \frac{n (n/e)^n}{((n+n_0)/e)^{n+n_0}}\le \frac{e^{n_0}}{n^{n_0-1}}
\end{equation*}
for $n,n_0\in\mathbb N$. Now fix ${\bs t}\in \mathbb R^M$. We estimate
\begin{equation*}
\begin{split}
\left|\frac{m_{\bs k}{\bs t}^{\bs k}}{{\bs k}!} \right|&=
\frac{\Gamma(\gamma_0)}{\Gamma(\gamma_{\bs k})}
\frac{|f_{{\bs k},0}|}{f_{{\bf 0},0}u_c^{\gamma_{\bs k}-\gamma_{\bf 0}}} |{\bs t}^{\bs k}|
\le
\frac{\Gamma(\gamma_0)}{(|{\bs k}|+\lfloor \frac{|{\bs k}|-3}{2}\rfloor)!}
\frac{D |{\bs k}|! {\bs R}^{\bs k}}{f_{{\bf 0},0}} |{\bs t}^{\bs k}|
\max (u_c,u_c^{-1})^{\gamma_{\bs k}-\gamma_{\bf 0}}\\
&\le \frac{D\Gamma(\gamma_0)}{f_{{\bf 0},0}} e^{|\bs k|}{\bs R}^{\bs k}
|{\bs t}^{\bs k}|\frac{1}{|{\bs k}|^{\lfloor \frac{|{\bs k}|-5}{2}\rfloor}}
\max (u_c,u_c^{-1})^{\frac{M+2}{2}|{\bs k}|}.
\end{split}
\end{equation*}
The r.h.s. vanishes as $|{\bs k}|\to\infty$, which implies the assertion.
The equivalent statement is obvious.
\end{proof}

Our proof of claim \textit{i)} in Theorem~\ref{theo:probdist} relies
on an application of L\'evy's continuity theorem \cite[Thm.~23.8]{B96},
which we cite in the following lemma.

\begin{lemma}[L\'evy's continuity theorem \cite{B96}]\label{lem:Levy}
For $n\in\mathbb N$, let probability measures $\mu_n$ on the Borel $\sigma$-algebra 
of $\mathbb R^M$ be given. If the sequence $\{\widehat \mu_n\}_{n\in\mathbb N}$ of 
their characteristic functions $\widehat \mu_n$ converges pointwise to a complex 
function $\phi$ which is continuous at the origin, then $\phi$ is the characteristic 
function of a uniquely determined Borel probability measure $\mu$. Moreover, 
$\{\mu_n\}_{n\in\mathbb N}$ converges to $\mu$ weakly. \qed 
\end{lemma}

\begin{proof}[Proof of claim \textit{i)} in Theorem~\ref{theo:probdist}.]
Lemma \ref{lem:pos} implies that the Borel probability measures
$\mu_{n_0}$ eqn.~\eqref{form:normrv} are well-defined for almost all 
$n_0\in\mathbb N$. The estimate in Lemma \ref{lem:momcon} implies uniform 
convergence of the sequence of Fourier transforms $\widehat \mu_{n_0}:\mathbb 
R^M\to\mathbb C$ of $\mu_{n_0}$, in every ball of finite radius centred 
at the origin. In particular, we have pointwise convergence of the sequence 
$\{\widehat \mu_{n_0}\}_{n_0\in\mathbb N}$. For $M=1$, 
the corresponding argument is given in the proof of \cite[Thm~6.4.5]{C74}. 
It can be directly extended to arbitrary $M$. Since the functions 
$\widehat \mu_{n_0}$ are continuous, we conclude that the limit function 
$\phi:\mathbb R^M\to\mathbb C$ is continuous at the origin.
Now apply L\'evy's continuity theorem. The limit probability measure $\mu$ 
has moments $m_{\bs k}$ eqn.~\eqref{form:normmom}. The claimed statement 
of the theorem follows, when phrasing the result in terms of the associated 
random variables.
\end{proof}

The connection to Brownian motion, which is claimed in part \textit{ii)} of 
Theorem~\ref{theo:probdist}, will be established in the following section.

\section{Dyck paths revisited}\label{sec:Dyckrev}

We apply our results to the example of Dyck paths of Section \ref{sec:Dyck}. 
The power series solution $E({\bs u_0})$ of eqn.~\eqref{form:exfunc}, specialised to 
$\bs u=\bs u_0$, has radius of convergence $u_c=1/4$, with a square root 
singularity at $u=u_c$, and $E({\bs u_c})=1$. The $q$-functional equation 
eqn.~\eqref{form:exfunc} satisfies Assumption \ref{ass}. Since the random variables 
$(X_{1,n_0},\ldots,X_{M,n_0})$ of eqn.~\eqref{form:DyckRV} have the distribution of 
$\mu_{n_0}$, as defined in eqn.~\eqref{form:normrv}, Proposition 
\ref{prop:dyckmom} follows from the results of the previous section.

\begin{proof}[Proof of Proposition \ref{prop:dyckmom}]
The generating function $E({\bs u})$ of Dyck paths satisfies the 
$q$-functional equation eqn.~\eqref{form:exfunc}. Assumption \ref{ass} holds with 
$u_c=1/4$ and $y_c=1$. We have $f_{{\bf 0},0}=-4$, $\gamma_{\bf 0}=-1/2$, 
$\mu_i=(i+1)/4$ for $i=1,\ldots,M-1$ and $\mu_0=1/8$. Thus, Proposition 
\ref{prop:rek} yields eqn.~\eqref{form:exc2}. The distribution of the 
random variables $(X_{1,n_0},\ldots,X_{M,n_0})$ of eqn.~\eqref{form:DyckRV} is that of 
the probability measures $\mu_{n_0}$ in eqn.~\eqref{form:normrv}. By 
part \textit{i)} of Theorem \ref{theo:probdist}, there exists a unique limit probability 
measure $\mu$. We thus get eqn.~\eqref{eqn:dyckres} from 
eqn.~\eqref{form:normmom}.
\end{proof}

The connection between Dyck paths and Brownian excursions in 
Proposition~\ref{prop:reldyckexc} leads, for a general 
$q$-functional equation, to an explicit characterisation of the 
limit probability measure $\mu$, resp.~of the associated limit random variable 
$(Y_1,\ldots,Y_M)$.

\begin{proof}[Proof of claim \textit{ii)} in Theorem~\ref{theo:probdist}.]
For the given $q$-functional equation, let $F_0(\bs\epsilon)$ be the 
generating function of the leading amplitudes $f_{\bs k,0}$. For $k=1,\ldots,M$ 
and $d_k\in\mathbb R$, define 
$G_0(\bs\epsilon)=F_0(\epsilon_1 d_1,\ldots, \epsilon_M d_M)$.
An explicit calculation using eqn.~\eqref{form:ee1} shows that the power series 
$G_0(\bs\epsilon)$ satisfies the same type of differential equation 
as $F_0(\bs\epsilon)$ does, with $\mu_0$ replaced by $\mu_0 d_1$ and 
$\mu_i$ replaced by $\mu_i d_{i+1}/d_i$, where $i=1,\ldots,M-1$. Now choose 
the values $d_k$, such that the equation for $G_0(\bs\epsilon)$ is 
that of Dyck paths. Noting the relation between Dyck paths and Brownian 
excursions eqn.~\eqref{form:DyEx}, we arrive at the values $c_k=2^{(k+2)/2}
/d_k$, for numbers $c_k$ as in the claim of Theorem~\ref{theo:probdist}.
\end{proof}

\section{Concluding remarks}\label{sec:con}

We finally stress three central aspects of our approach.
Firstly, the approach yields a \emph{universal} limit distribution -- loosely
spoken, it appears for all models, whose underlying functional equation 
has the same singularity structure. Such a result may be compared to a 
central limit theorem in probability 
theory. For example, models other than Dyck paths, which display a square root 
as dominant singularity of their size generating function, are certain models of 
trees or polygons. For simply generated trees, $q$-functional equations appear when counting 
by number of vertices and by internal path length \cite{T91}. Using the 
above setup, moment recursions for the parameter ``sum of $k$-th 
powers of the vertex distances to the root'' are obtained. For polygon 
models, $q$-functional equations appear when counting by perimeter 
and area \cite{D99,R02}, which is, for column-convex polygons, the 
sum of the column heights. The above setup gives moment recursions
for the parameter ``sums of $k$-th powers of the column heights'', in the 
limit of infinite (horizontal) perimeter.

Secondly, our approach is \emph{algorithmic} -- the moment recursion can
finally be deduced from a straightforward calculation, by the method of
dominant balance. Our approach also allows to obtain corrections to the 
asymptotic behaviour, which cannot (easily) be deduced by other methods.

Thirdly, the approach is \emph{flexible} -- it may be applied to other classes 
to obtain a universal limit distribution, which only depends on the singularity 
structure of the functional equation. Our generating function approach is particularly
suited for counting parameters, which decompose linearly under the 
cartesian product construction. Examples of such models with a 
rational generating function appear in \cite{R02}. Examples with an inverse 
square root appear in \cite{NT03}. In particular, the discrete counterparts of 
Brownian motion, Brownian bridges, and Brownian meanders can be studied, 
see \cite{NT03,R06}. Also, universality questions for parameters related 
to left and right path lengths in trees \cite{J06,KS06,P06,Bou06} can be studied 
by our methods, compare the discussion in \cite{BJ06}.

Within the framework of simply generated trees, an alternative derivation of 
the above moment recursion could be obtained with the techniques of 
\cite{J03}, where the different problem of the (generalised) Wiener index of 
trees was analysed. It would be interesting to 
consider how our methods can be adapted to this problem.

\section*{Acknowledgements}

The author thanks Philippe Duchon, Philippe Flajolet and Michel 
Nguy$\tilde{\mbox{\rm \^e}}$n Th$\acute{\mbox{\rm \^e}}$ for helpful 
discussions, and Svante Janson for comments on the manuscript. 
He thanks the department LaBRI (Bordeaux) for hospitality 
in autumn 2003, where parts of the problem have been analysed. The
author thanks the referees for suggestions, which improved the presentation 
of the article. Financial support by the German Research Council (DFG) 
is gratefully acknowledged.

\end{document}